\documentclass{amsproc}

\usepackage[all]{xy}

\usepackage{hhline}
\usepackage{amssymb}
\usepackage[dvips, dvipsnames, usenames]{color}
\usepackage{eucal}

\newcommand{\PSL}{\mathbf{PSL}}

\newcommand{\ydf}{{}_{\ku F}^{\ku F}\mathcal{YD}}
\newcommand{\ydgx}{{}_{\ku \Gx}^{\ku \Gx}\mathcal{YD}}
\def\ixq{\Inn_{X,\q}}

\newcommand{\X}{\tt{X}}

\newcommand\sn{\mathbb S_n}
\newcommand\bn{\mathbb B_n}

\newcommand{\fd}{finite-dimensional}
\newcommand{\diag}{\operatorname{diag}}

\newcommand{\Fun}{\operatorname{Fun}}

\newcommand{\Imm}{\operatorname{Im}}

\newcommand{\ord}{\operatorname{ord}}

\newcommand{\Inn}{\operatorname{Inn}}

\def\gax{\Inn_{\trid} (X)}

\def\Gx{\mathbb G_X}

\newcommand{\cl}{{\mathcal C}_{\ell}}

\newcommand{\fp}{{\mathbb F}_{\hspace{-2pt}p}}
\newcommand{\kc}{\mathbb F_q}

\renewcommand{\_}[1]{\mbox{$_{\left( #1 \right)}$}}

\newcommand\toba{{\mathfrak B }}

\newcommand{\gr}{\operatorname{gr}}
\newcommand{\trid}{\triangleright}

\newcommand{\Lb}{{\mathbb L}}

\newcommand{\R}{{\mathcal R}}

\newcommand{\ku}{\mathbb C}

\newcommand{\Kb}{{\mathbb K}}

\newcommand{\Z}{{\mathbb Z}}
\newcommand{\N}{{\mathbb N}}

\newcommand{\G}{{\mathbb G}}

\newcommand{\Q}{{\mathsf Q}}

\newcommand{\C}{{\mathcal C}}
\newcommand{\D}{{\mathcal D}}

\newcommand{\q}{{\mathbf q}}

\newcommand{\GL}{\mathbf{GL}}
\newcommand{\PGL}{\mathbf{PGL}}
\newcommand{\SL}{\mathbf{SL}}

\newcommand{\T}{{\mathcal T}}
\newcommand{\Hc}{{\mathcal H}}

\newcommand{\Oc}{{\mathcal O}}
\newcommand{\oc}{{\mathcal O}}

\newcommand{\ydh}{{}^{\Hc}_{\Hc}\mathcal{YD}}

\newcommand{\ydhs}{{}^{\Hc_{\phi}}_{\Hc_{\phi}}\mathcal{YD}}

\newcommand{\ydg}{{}^{\ku G}_{\ku G}\mathcal{YD}}

\newcommand{\Ss}{{\mathcal S}}

\newcommand{\Aut}{\operatorname{Aut}}
\newcommand{\Out}{\operatorname{Out}}

\newcommand\card{\operatorname{card}}

\theoremstyle{plain}

\newtheorem{lema}{Lemma}[section]
\newtheorem{theorem}[lema]{Theorem}
\newtheorem{cor}[lema]{Corollary}

\newtheorem{prop}[lema]{Proposition}

\newtheorem{question}{Question}
\newtheorem{question-app}{Question}

\newtheorem*{ques}{Question}

\newtheorem*{strategy}{Strategy}
\newtheorem*{alg}{Algorithm}

\theoremstyle{definition}
\newtheorem{definition}[lema]{Definition}
\newtheorem{exa}[lema]{Example}

\theoremstyle{remark}
\newtheorem{obs}[lema]{Remark}

\newcommand\id{\operatorname{id}}

\newcommand\st{\mathbb S_3}
\newcommand\sk{\mathbb S_4}
\newcommand\sco{\mathbb S_5}

\newcommand\an{\mathbb A_n}

\newcommand\am{\mathbb A_m}

\newcommand\ac{\mathbb A_4}
\newcommand\aco{\mathbb A_5}
\newcommand\as{\mathbb A_6}
\newcommand\A{\mathbb A}

\newcommand\Sim{\mathbb S}

\newcommand\sm{\mathbb S_m}

\newcommand\s{\mathbb S}

\def\pf{\begin{proof}}
\def\epf{\end{proof}}

\theoremstyle{remark}

\begin{document}

\renewcommand{\baselinestretch}{1.2}

\thispagestyle{empty}

\title[On Nichols algebras associated to simple racks]
{On Nichols algebras associated to simple racks}

\author[Andruskiewitsch, Fantino, Garc\'\i a, Vendramin]{N. Andruskiewitsch,
F. Fantino, G. A. Garc\'\i a, L. Vendramin}

\address{\noindent N. A., F. F., G. A. G.:
Facultad de Matem\'atica, Astronom\'{\i}a y F\'{\i}sica,
Universidad Nacional de C\'ordoba. CIEM -- CONICET. 
Medina Allende s/n (5000), Ciudad Universitaria, C\'ordoba,
Argentina
\newline
\indent F. F., G. A. G.:
Facultad de Ciencias Exactas, F\'isicas y Naturales,
Universidad Nacional de C\'ordoba. 
Velez Sarsfield 1611 (5000), Ciudad Universitaria, C\'ordoba,
Argentina
\newline \indent L. V. : Departamento de Matem\'atica, FCEyN,
Universidad de Buenos Aires, Pab. I,  Ciudad Universitaria (1428),
Buenos Aires, Argentina
\newline  \indent L. V. : Instituto de Ciencias,
Universidad de Gral. Sarmiento, J.M. Gutierrez
1150, Los Polvorines (1653), Buenos Aires, Argentina}

\email{(andrus, fantino, ggarcia)@famaf.unc.edu.ar}
\email{lvendramin@dm.uba.ar}

\thanks{This work was partially supported by ANPCyT-Foncyt,
CONICET, Ministerio de Ciencia y
Tecnolog\'{\i}a (C\'ordoba), Secyt-UNC and Secyt-UBA}

\date{\today}

\dedicatory{Dedicado a Mat\'\i as Gra\~na}

\begin{abstract}
This is a report on the present state of the problem of
determining the dimension of the Nichols algebra associated to a
rack and a cocycle. This is relevant for the classification of
finite-dimensional complex pointed Hopf algebras whose group of
group-likes is non-abelian. We deal mainly with simple racks. We
recall the notion of rack of type D, collect the known lists of
simple racks of type D and include preliminary results for the
open cases. This notion is important because the Nichols algebra
associated to a rack of type D and any cocycle has infinite
dimension. For those racks not of type D, the computation of the
cohomology groups is needed. We discuss some techniques for this
problem and compute explicitly the cohomology groups corresponding
to some conjugacy classes in symmetric or alternating groups of
low order.
\end{abstract}
\maketitle

\section{Introduction}
Throughout the paper we work over the field $\ku$ of complex numbers.
The problem of classifying finite-dimensional pointed Hopf
algebras over non-abelian finite groups reduces in many cases to a
question on conjugacy classes. In this introduction we give a
historical account and place the problem in the overall picture.

\medbreak\subsection{} We briefly recall the lifting method for the
classification of pointed Hopf algebras, see
Subsection \ref{subsec:nichols} for unexplained terminology and
\cite{AS-cambr} for a full exposition.
Let $H$ be a Hopf
algebra with bijective antipode and assume that the coradical
$\displaystyle H_0 = \sum_{C \text{ simple subcoalgebra of } H} C$
is a Hopf subalgebra of $H$.
Consider the coradical filtration of $H$:
$$
H_0 \subset H_1 \subset \dots \subset H = \bigcup_{n\geq 0} H_n,$$
where $H_{i+1} = \{x\in H: \Delta(x) \in H_i\otimes H + H\otimes H_0\}$.
Then the associated graded coalgebra
$\gr H$ has a decomposition $\gr H\simeq R \# H_0$, where $R$ is an algebra
with some special properties and
$\#$ stands for a kind of semidirect product (technically,
a Radford biproduct or bosonization; the underlying vector space is $R \otimes H_0$).
The algebra $R$, more precisely, is a Hopf algebra in the braided
tensor category of Yetter-Drinfeld modules over $H_0$,
see Subsection \ref{subsec:nichols}, and inherits the grading
of $\gr H$: $R = \oplus_{n\geq 0} R^n$. If $V = R^1$, then
the subalgebra of $R$ generated by $V$ is isomorphic to the
Nichols algebra $\toba(V)$, that is completely determined by
the Yetter-Drinfeld module $V$.

Let us fix a semisimple Hopf algebra $A$. One of the fundamental
steps of the lifting method to classify \fd {} Hopf algebras $H$ with
$H_0 \simeq A$ is to address the following question, see \cite{A}:

\begin{ques}
Determine the Yetter-Drinfeld modules $V$ over $A$ such that the
dimension of $\toba(V)$ is finite, and if so, give an efficient
set of relations of $\toba(V)$.
\end{ques}

An important observation is that the Nichols algebra $\toba(V)$,
as algebra and coalgebra, is completely determined just by
the braiding $c: V\otimes V \to V \otimes V$. Therefore, it is
convenient to consider classes of braided vector spaces $(V, c)$
depending on the class of semisimple Hopf algebras we are considering.

\medbreak\subsection{} A Hopf algebra $H$ is pointed if $H_0$ is
isomorphic to the group algebra $\ku G$,
where $G$ is the group of grouplikes of $H$. Let us consider
first the case when $G$ is abelian.
A braided vector space $(V,c)$ is of diagonal type if $V$ has a basis
$(v_i)_{1\le i \le n}$ such that $c(v_i \otimes v_j) =q_{ij}
v_j \otimes v_i$, where the $q_{ij}$'s are non-zero scalars \cite{AS-adv}.
The class of braided vector spaces of diagonal type
corresponds to the class of pointed Hopf algebras
with $G$ abelian (and finite). A remarkable result is the complete list of all
braided vector spaces of diagonal type with
finite-dimensional Nichols algebra  \cite{H-all};
the basic tool in the proof of this result is the
Weyl groupoid \cite{H-weyl}.
The classification of all \fd{}  pointed Hopf algebras
with $G$ abelian and order of $G$ coprime with 210
was obtained in \cite{AS-05}, relying crucially on \cite{AS-adv, H-all}.
Notice however that the article \cite{H-all} does
not contain the efficient set of relations for
finite-dimensional Nichols algebras of diagonal type;
so far, this is available for the special classes of braided
vector spaces of Cartan type \cite{AS-adv} and more
generally of standard type \cite{Ang}.

\medbreak\subsection{} Let us now turn to the case when
$H$ is pointed  with $G$ non-abelian and mention some antecedents.

\begin{itemize}
\renewcommand{\labelitemi}{$\diamond$}

\medbreak  \item The first genuine examples of \fd{}
pointed Hopf algebras with non-abelian group appeared in \cite{MS, FK},
  as bosonizations of  Nichols algebras related to the transpositions
in $\st$ and $\sk$, see Subsection \ref{subsec:fk}.
  The analogous quadratic algebra over $\sco$ was computed by
Roos with a computer and proved to be a Nichols algebra in \cite{G-zoo}.

\medbreak  \item In  \cite{G-cm},
Gra\~na identified the class of braided vector
spaces corresponding to pointed Hopf algebras
with non-abelian group as those constructed from racks and
cocycles. He also computed in \cite{G-zoo}
several \fd{} Nichols algebras
with the help of computer programs.

\medbreak  \item In \cite{G-cm}, Gra\~na also suggested to
look at braided vector subspaces to decide  that a Nichols
algebra has infinite dimension. After \cite{H-all}, this idea
was implemented in several papers, by looking at abelian
subracks. See \cite{AF1, AF2, afz, az, fantino-2007, FGV, FGV2}.

\medbreak  \item The construction of the Weyl groupoid for braided vector
spaces of diagonal type in \cite{H-weyl} was extended to
braided vector spaces arising from semisimple Yetter-Drinfeld modules
in \cite{AHS}. This allowed to consider braided vector subspaces associated to non-abelian subracks \cite{AF3}. 
A further study of the Weyl groupoid in \cite{AHS} was undertaken in \cite{HS1}.
An important consequence of one of the results in \cite{HS1}
is the notion of rack of type D \cite{AFGV}.
\end{itemize}

\medbreak\subsection{} We shall explain in detail the notion of
rack of type D in Subsection \ref{subsec:racks-typeD}, but we try
now to give a glimpse. As we explain in Subsection
\ref{subsec:nichols}, our goal is to determine if the Nichols
algebra $\toba(\oc, \rho)$ related to a conjugacy class $\oc$ in a
finite group $G$ and a representation $\rho$ of the centralizer is
\fd{}. We say that the conjugacy class $\Oc$ is \emph{of type D}
if there exist $r,s\in \Oc$ such that
\begin{enumerate}
  \item  $(rs)^2 \neq (sr)^2$,
  \item  $r$ is not conjugated to $s$ in the subgroup of $G$ generated by $r,s$.
\end{enumerate}

Then $\dim \toba(\oc, \rho) = \infty$ \emph{for any} $\rho$; furthermore
this will happen for any group
$G'$ containing $\Oc$ as a conjugacy class (that is, as a subrack).
By reasons exposed in Subsection \ref{subsec:racks-typeD},
we focus on the following case.

\begin{question}\label{que:rackssimples-typeD}
Determine all simple racks of type D.
\end{question}

The classification of finite simple racks is known, see
Subsection \ref{subsec:simple-racks};
the list consists
of conjugacy classes in groups of 3 types. In other words,
we need to check, for each conjugacy class in the list of simple racks,
whether there exist $r,s$ satisfying (1) and (2) above.
The main purpose of this paper is to report the actual
status of this purely group-theoretical question,
that is succinctly as follows.

\begin{itemize}
\renewcommand{\labelitemi}{$\diamond$}

\medbreak  \item \cite{AFGV} The conjugacy classes in the
alternating and symmetric groups, $\am$ and $\sm$, are of type D, except for
a short list of exceptions listed in
Theorems \ref{th:racks-an-liquidados} and \ref{th:racks-sn-liquidados};
for some of these exceptions, we know
that they are not of type D, see Remark 4.2 in \emph{loc. cit}.

\medbreak  \item \cite{AFGV-espo}
 The conjugacy classes in the sporadic groups are of type D, except for
a short list of exceptions listed in Theorems \ref{th:racks-liquidados};
for some of these, we know that they are not of type D, see Table \ref{tab:0}.
The verification was done with the help of  \textsf{GAP}, see \cite{logbook}.

\medbreak  \item \cite{FV} Twisted conjugacy
classes of sporadic groups are also mostly of type D,
except for a short list of exceptions, see
 Theorem \ref{th:esporadicos-twisted-liquidados}.

\medbreak  \item \cite{AFGV-thr} Some techniques
to deal with twisted homogenous racks were found; so far,
most of the examples dealt with  are of type D.

\medbreak  \item We include in Subsection
\ref{subsect:-finite-lie} some preliminary
results on conjugacy classes on simple groups of Lie type;
again, most of the  examples are of type D.

\medbreak  \item The simple affine racks do not seem to be of type D.
\end{itemize}

What happens beyond type D? As we see by now,
there are roughly two large classes of simple racks,
one formed by the affine ones and the conjugacy
class of transpositions in $\sm$, and the rest. For this second class,
our project is to finish the determination of those
of type D and attack the remaining ones as explained on page
\pageref{strategy:racks}. That is, to compute the
pointed sets of cocycles of degree $n$ and then try to discard
the corresponding braided vector spaces by abelian techniques.
The first class is not tractable by the strategy of subracks.
We should also mention the recent paper \cite{GHV}
with a different approach.

\medbreak\subsection{} The paper is organized as follows. We
discuss Nichols algebras, racks, cocycles, the criterion of type
D, the classification of finite simple racks and the strategy of
subracks in Section \ref{sec:preliminaries}. Section
\ref{sec:tools} contains some techniques for the computation of
cocycles. In the next sections we list explicitly the simple racks
that are known to be of type D. In Section
\ref{sec:applications-pointed} we illustrate the consequences of
these results to the classification of pointed Hopf algebras. In
Appendix \ref{subsec:nichols-examples-fd}, we list all known
examples of \fd{} Nichols algebras associated to racks and
cocycles; in Appendix \ref{subsec:questions}, we put together some
questions scattered along the text. 

This survey contains also a
few new concepts and results, among them: the computation of the
enveloping group of the rack of transpositions in $\sm$, see
Proposition \ref{prop:enveloping-group-sm}; the twisting operation
for cocycles on racks, see Subsection \ref{subsec:twisting}; the
calculation of some cohomology groups using the program
\textsf{RiG}, see Subsection \ref{subsec:rig}; some preliminary
discussions on conjugacy classes of type D in finite groups of Lie
type, see Subsection \ref{subsect:-finite-lie}.

\section{Preliminaries}\label{sec:preliminaries}

\subsection*{Conventions} $\N= \{1, 2, 3, \dots\}$; $\s_X := \{f: X \to X$ bijective$\}$; if $m\in \N$, then $\G_m$ is the group of $m$-th roots of 1 in $\ku$.

\subsection{Racks}\label{subsect:racks}

\

 We briefly recall the basics of racks; see \cite{AG-adv}
for more information and references. A \emph{rack} is a pair
$(X,\trid)$ where $X$ is a non-empty set and $\trid:X\times X\to
X$ is an operation such that
\begin{eqnarray}
    &\text{the map }\varphi_x = x\trid \underline{\quad}
        \text{ is bijective for any }x\in X,&\ \text{and} \\
    &x \trid(y\trid z) = (x\trid y) \trid(x\trid z)
        \text{ for all }x,y,z\in X.&\label{eqn:selfdist}
\end{eqnarray}

A group $G$ is a rack with $x\trid y = xyx^{-1}$, $x,y\in G$; if
$X\subset G$ is stable under conjugation by $G$, that is a union
of conjugacy classes, then it is a subrack of $G$. The main idea
behind the consideration of racks is to keep track just of the
conjugation of a group. Morphisms of racks and subracks are
defined as usual. For instance, $\varphi:X\to\s_X$, $x\mapsto
\varphi_x$, is a morphism of racks, for any rack $X$. Any rack $X$
considered here satisfies the conditions
\begin{align}
\label{eqn:quandle} x\trid x &= x, \\
\label{eqn:crossed-set} x\trid y &= y \implies y\trid x = x,
\end{align}
for any $x, y\in X$. This is technically a \textit{crossed set}, but we
shall simply say a rack. So, we rule out, for example, the
permutation rack $(X, \sigma)$, where $\sigma \in \mathbb S_X$ and
$\varphi_x =\sigma$ for all $x$.

\medbreak  The rack with just one element is called \emph{trivial}.

\medbreak  We shall consider some special classes
of racks that we describe now.

\medbreak \emph{Affine racks}. If $A$ is an abelian group and
$T\in\Aut(A)$, then $A$ is a rack with $x\trid y=(1-T)x+Ty$. This
is called an \textit{affine rack} and denoted $\Q_{A,T}$.

\medbreak \emph{Twisted conjugacy classes}. Let $G$ be a finite
group and $u \in \Aut (G)$; $G$ acts  on itself by
$x\rightharpoonup_u y = x\,y\,u(x^{-1})$, $x,y\in G$.
The orbit $\oc^{G,u}_x$ of $x\in G$ by this action is a rack
with operation
\begin{equation}\label{eqn:twisted-conjugacy-class}
y\trid_u z = y\,u(z\,y^{-1}), \quad y,z\in \oc^{G,u}_x.
\end{equation}
We shall say that $\oc^{G,u}_x$ is a \textit{twisted conjugacy class} of type
$(G,u)$.

\subsubsection*{Notation.}
\begin{itemize}    \item $\T = $ any of the conjugacy classes of 3-cycles in $\ac$
    (the tetrahedral rack).

\smallbreak\item $\Q_{A,T} =$ affine rack associated to an abelian group  $A$  and
$T\in\Aut(A)$.

\smallbreak\item $\D_n = $  class of involutions in the dihedral
    group of order $2n$, $n$ odd.

\smallbreak\item $\Oc^m_j = $ conjugacy class of $j$-cycles in
$\sm$.
\end{itemize}

\medbreak We need some terminology on racks.

\begin{itemize}
    \item A rack $X$ is \emph{decomposable} if it can expressed as a disjoint
    union of subracks $X=X_1\coprod X_2 $. Otherwise,  $X$ is \emph{indecomposable}.


    \medbreak    \item A rack $X$ is said to be \emph{simple} iff $\card X>1$ and for any surjective
        morphism of racks $\pi: X \to Y$, either $\pi$ is a bijection or $\card Y = 1$.

    \medbreak    \item If $X$ is a rack and $j\in \Z$, then
    $X^{[j]}$ is the rack with the same set $X$ and
    multiplication $\trid^j$ given by
$x\trid^j y = \varphi_x^{j}(y)$, $x,y\in X$.

\end{itemize}

\subsection{Nichols algebras}\label{subsec:nichols}

\

Nichols algebras play a crucial role in the classification of Hopf
algebras, see \cite{AS-cambr} or a brief account in
Section \ref{sec:applications-pointed} below.
Let $n \geq 2$ be an integer. We start by reminding the well-known
presentations by generators and relations of the braid group
$\bn$ and the symmetric group $\sn$:
\begin{align*}
\bn &= \langle (\sigma_i)_{1\le i \le n-1}\vert
\sigma_i\sigma_j\sigma_i = \sigma_j\sigma_i\sigma_j, \, \vert i-j\vert =1;
\quad \sigma_i\sigma_j = \sigma_j\sigma_i, \, \vert i-j\vert >1\rangle
\\ \sn &= \langle (s_i)_{1\le i \le n-1} \vert  s_is_js_i = s_js_is_j, \, \vert i-j\vert =1;
 \quad s_is_j = s_js_i, \, \vert i-j\vert >1; \quad  s_i^2 = e\rangle,
\end{align*}
indices in the relations going over all possible $i,j$. There is a canonical projection
$\pi:\mathbb B_n \to \sn$, that admits a so-called Matsumoto section
$M: \sn \to \bn$; this is not a morphism of groups, and
it is defined by $M(s_i) = \sigma_i$, $1\le i \le n-1$, and $M(st) = M(s)M(t)$,
for any $s,t\in \sn$ such that $l(st) = l(s) + l(t)$,
$l$ being the length of a word in generators $s_i$. 

Let $V$ be a vector space and $c\in \GL(V \otimes V)$. Recall that $c$ fulfills the braid
equation if $(c\otimes \id)(\id\otimes c)(c\otimes \id) = (\id\otimes c)(c\otimes \id)(\id\otimes c)$.
In this case, we say that $(V,c)$ is a \emph{braided vector space} and that $c$ is a \emph{braiding}.
Since $c$ satisfies the braid equation, it induces a representation of the braid group $\bn$,
 $\rho_n: \bn \to \GL(V^{\otimes n})$, for each $n\ge 2$. Explicitly, $\rho_n(\sigma_i)
= \id_{V^{\otimes (i-1)}} \otimes c \otimes \id_{V^{\otimes (n -i-1)}}$, $1\le i \le n-1$.
 Let
 \begin{equation}\label{eqn:quantum-symmetrizer}
 Q_n = \sum_{\sigma \in \sn} \rho_n(M(\sigma))\in End (V^{\otimes n}).
 \end{equation}
Then the \emph{Nichols algebra }$\toba(V)$ is the quotient of the
tensor algebra $T(V)$ by $\oplus_{n\ge 2} \ker Q_n$, in fact a
2-sided ideal of $T(V)$. If $c = \tau$ is the usual switch, then
$\toba(V)$ is just the symmetric algebra of $V$; if $c = -\tau$,
then $\toba(V)$ is the exterior algebra of $V$. But the
computation of the Nichols algebra of an arbitrary braided vector
space is a delicate issue. We are interested in the Nichols
algebras of the braided vector spaces arising from Yetter-Drinfeld
modules\footnote{Any braided vector space with \emph{rigid}
braiding arises as a Yetter-Drinfeld module \cite{Tk}.}.

A \emph{Yetter-Drinfeld module} over a Hopf algebra $H$ with
bijective antipode $\Ss$ is a left $H$-module $M$ and
simultaneously a left $H$-comodule, with coaction $\lambda: M \to
H \otimes M$ compatible with the action in the sense that $\lambda
(h\cdot x) = h\_1 x\_{-1} \Ss(h\_3) \otimes h\_2\cdot x\_{0}$, for
any $h\in H$, $x\in M$. Here $\lambda (x) = x\_{-1}  \otimes
x\_{0}$, in Heyneman-Sweedler notation. A Yetter-Drinfeld module
$M$ is a braided vector space with $c(m\otimes n) = m\_{-1}\cdot n
\otimes m\_0$, $m,n\in M$. We shall see in Section
\ref{sec:applications-pointed} how Nichols algebras of
Yetter-Drinfeld modules enter into the classification of Hopf
algebras. In this paper, we are interested in the case when $H=
\ku G$ is the group algebra of a finite group $G$. In this
setting, a Yetter-Drinfeld module over $H$ is a left $G$-module
$M$ that bears also a $G$-grading $M = \oplus_{g\in G}M_g$,
  compatibility meaning that $h\cdot M_g = M_{hgh^{-1}}$
for all $h,g\in G$; the braiding is
  $c(m\otimes n) = g\cdot n \otimes m$, $m\in M_g$, $n\in M$.

Now a braided vector space may be realized as a Yetter-Drinfeld
module over many different groups and in many different ways. It
is natural to look for a description of the class of braided
vector spaces that actually arise as Yetter-Drinfeld modules over
some finite group and to study them by their own. If $G$ is a
finite group, then any Yetter-Drinfeld module over the group
algebra $\ku G$ is semisimple. Furthermore, it is well-known that the
set of isomorphism classes of irreducible Yetter-Drinfeld modules
over $\ku G$ is parameterized by pairs $(\Oc, \rho)$, where $\Oc$
is a conjugacy class of $G$ and $\rho$ is an irreducible
representation of the centralizer of a fixed point in $\Oc$. M.
Gra\~na observed that the class of braided vector spaces arising
from Yetter-Drinfeld modules over finite groups is described using
racks and cocycles, see \cite{G-cm} and also \cite[Th.
4.14]{AG-adv}.

\subsection{Nichols algebras associated to racks and cocycles}\label{subsec:nichols-racks}

\

We are focused in this paper on
Nichols algebras associated to braided vector spaces built from
racks and cocycles. We start by describing the cocycles associated to racks.
Let $X$ be a rack and $n\in \N$. A map $\q:X\times X\to\GL(n,\ku)$
is a \emph{2-cocycle} of degree $n$ if
$$\q_{x,y\trid z}\q_{y,z}= \q_{x\trid y,x\trid z}\q_{x,z},$$ for all
$x,y,z\in X$. Let $\q$ be a 2-cocycle of degree $n$, $V= \ku X\otimes\ku^{n}$,
where $\ku X$ is the vector space with basis $e_x$, for $x\in X$.
We denote $e_xv := e_x\otimes v$. Consider the linear isomorphism
$c^{\q}:V\otimes V\to V\otimes V$,
\begin{equation}\label{eqn:trenza}
c^{\q}(e_xv\otimes e_yw)=e_{x\trid y}\q_{x,y}(w)\otimes e_xv,
\end{equation}
$x,y \in X$, $v,w\in\ku^{n}$. Then $c^{\q}$ is
a solution of the braid equation:
\begin{align*}
(c^{\q}\otimes \id)(\id\otimes c^{\q})(c^{\q}\otimes \id) =
(\id\otimes c^{\q})(c^{\q}\otimes \id)(\id\otimes c^{\q}).
\end{align*}

\begin{exa}\label{exa:inverse-rack}
Let $X$ be a finite rack and $\q$ a 2-cocycle.
The dual braided vector space of $(\ku X\otimes\ku^{n}, c^{\q})$
is isomorphic to $(\ku X^{[-1]}\otimes\ku^{n}, c^{\widehat{\q}})$
where $\widehat{\q}_{x,y} = \q_{x, x\trid^{-1}y}$, $x, y \in X^{[-1]}$.
See Subsection \ref{subsect:racks} for $X^{[-1]}$.
\end{exa}

\medbreak The Nichols algebra associated to $c^{\q}$ is denoted
$\toba(X, \q)$.

\medbreak We need to consider only 2-cocycles (or simply cocycles,
for short) with some specific properties.

\begin{itemize}
    \item A cocycle $\q$ is \emph{finite} if its image generates a
finite subgroup of $\GL(n,\ku)$.

    \medbreak\item A cocycle $\q$ is \emph{faithful} if the morphism of
racks $g: X \to \GL(V)$ defined by $g_x(e_yw)= e_{x\trid
y} \q_{x, y}(w)$, $x,y\in X$, $w\in V$, is injective.
\end{itemize}

\medbreak We denote by $Z^{2}(X, \GL(n,\ku))$ the set of all
finite faithful 2-cocycles of degree $n$. Let $\q\in Z^{2}(X,
\GL(n,\ku))$ and $\gamma: X\to \GL(n,\ku)$ a map whose image generates a
finite subgroup. Define $\widetilde{\q}:X\times X\to\GL(n,\ku)$
\begin{equation}\label{eqn:rel-cociclos}
\widetilde{\q}_{ij} = \left(\gamma_{i\trid
j}\right)^{-1}\q_{ij}\gamma_{j}.
\end{equation}
Then $\widetilde{\q}$ is also a finite faithful cocycle and ``$\q
\sim \widetilde{\q}$ iff they are related by
\eqref{eqn:rel-cociclos} for some $\gamma$'' defines an
equivalence relation. We set
\begin{equation}\label{eq:ddcnc}
H^{2}(X, \GL(n,\ku)) = Z^{2}(X, \GL(n,\ku))/\sim.
\end{equation}

If $\q \sim \widetilde{\q}$, then the Nichols algebras $\toba(X,
\q)$ and $\toba(X, \widetilde{\q})$ are isomorphic as braided Hopf
algebras in the sense of \cite{Tk}, see \cite[Th. 4.14]{AG-adv}.
The converse is not true, see \cite{G-cm}.

\medbreak The main question we want to solve is the following.

\begin{question}\label{que:racks}
For any finite indecomposable rack $X$, for any $n\in \N$, and for
any $\q\in H^{2}(X, \GL(n,\ku))$, determine if $\dim \toba(X, \q)
< \infty$.
\end{question}

\begin{definition}\label{def:rack-collapses}
An indecomposable finite rack $X$ \emph{collapses at $n$} if for
any finite faithful cocycle $\q$ of degree $n$, $\dim\toba(X,\q)
=\infty$; $X$ \emph{collapses} if it collapses at $n$ for any
$n\in \N$.
\end{definition}

The first idea that comes to the mind is one would need to compute the group
$H^{2}(X, \GL(n,\ku))$ for any $n$. We shall see that in many
cases this is actually not necessary.

\begin{question}\label{que:collapses-degree1}
If $X$ collapses at $1$, does necessarily $X$ collapse?
\end{question}

Even partial answers to Question \ref{que:collapses-degree1} would
be very interesting.

\subsection{Racks of type D}\label{subsec:racks-typeD}

\

We now turn to a setting where the calculation of the cocycles is not needed.

\begin{definition}\label{def:tipod}
    Let $(X, \trid)$ be a rack. We say that $X$ is \emph{of type D} if there exists  a
    decomposable subrack $Y = R\coprod S$ of $X$ such that
    \begin{equation}\label{eqn:hypothesis-subrack}
        r\trid(s\trid(r\trid s)) \neq s, \quad \text{for some } r\in R, s\in S.
    \end{equation}
\end{definition}

The following important result is a consequence of \cite[Th. 8.6]{HS1},
proved using the main result of \cite{AHS}.

\begin{theorem}\label{th:racks-claseD} \cite[Th. 3.6]{AFGV}
    If $X$ is a finite rack of type D, then $X$  collapses. \qed
\end{theorem}

Therefore, it is very important to determine all simple racks of type D,
formally stated as Question \ref{que-app:rackssimples-typeD}.
The classification of simple racks is known and will be evoked below.
We focus on simple racks because of the following reasons:
    \begin{itemize}
\item If $Z$ is a finite rack and admits a rack epimorphism $\pi: Z\to X$,
where $X$ is of type D, then $Z$
            is of type D.

\medbreak\item If $Z$ is indecomposable, then it admits a rack
epimorphism $\pi: Z \to X$ with $X$ simple.
   \end{itemize}

\medbreak
We collect some criteria on racks of type D, see \cite[Subsection 3.2]{AFGV}.

    \begin{itemize}
\item If $Y\subseteq X$ is a subrack of type D, then $X$
is of type D.

\item If $X$ is of type D and $Z$ is a rack, then $X \times Z$ is
of type D.

\item  Let $K$ be a subgroup of a finite group $G$ and
 $\kappa\in C_G(K)$. Let $\R_{\kappa}: K\to G$ be the map given by $g\mapsto \widetilde g := g\kappa$.
    Let $\oc$, resp. $\widetilde{\oc}$,  be the conjugacy class of $\tau\in K$, resp. of $\widetilde{\tau}$ in $G$.
    Then $\R_{\kappa}$ identifies $\oc$ with a subrack of $\widetilde{\oc}$. Hence, if $\oc$ is of type D,
    then $\widetilde{\oc}$ is of type D.

   \end{itemize}

\medbreak There is a variation of the last criterium that needs the notion of \emph{quasi-real} conjugacy class.
Let $G$ be a finite group, $g \in G$ and $j \in \N$. Recall that $\oc^{G}_{g}$ is quasi-real of type $j$
if $g^{j}\neq g$ and $g^{j} \in \oc_{g}^{G}$.
If $g$ is real, that is $g^{-1} \in \oc_{g}^{G}$, but not an involution, then $\oc_{g}^{G}$ is quasi-real
of type $\ord(g) - 1$.

\begin{prop}\label{exa:jordan} \cite[Ex. 3.8]{AFGV}
Let $G$ be a finite group and $g = \tau\kappa\in G$, where $\tau\neq e$ and $\kappa\neq e$ commute.
Let $K =  C_G(\kappa) \ni \tau$; then
$\kappa\in C_G(K)$. Hence, the conjugacy class
$\oc$ of $\tau$ in $K$ can be identified with a subrack of the conjugacy
class $\widetilde{\oc}$ of $g$ in $G$ via $\R_\kappa$ as above.
Assume that
\begin{enumerate}
\item $\widetilde{\oc}$ and $\oc$ are quasi-real of type $j$,
\item the orders $N$ of $\tau$ and $M$ of $\kappa$ are coprime,
\item $M$ does not divide $j-1$,
\item there exist $r_0$, $s_0\in \oc$ such that $r_0\trid(s_0\trid(r_0\trid s_0)) \neq s_0$.
\end{enumerate}
Then $\widetilde{\oc}$ is of type D. \qed
\end{prop}

\subsection{Simple racks}\label{subsec:simple-racks}

\

Finite simple racks have been classified in \cite[Th. 3.9, Th.
3.12]{AG-adv}-- see also \cite{jo}. Explicitly, any simple rack
falls into one and only one of the following classes:

\medbreak
\begin{enumerate}

    \item \emph{Simple affine racks} $({\fp}^t, T)$, where $p$ a prime, $t\in \N$,  and
    $T$ is the companion matrix of a monic irreducible polynomial $f\in \fp[\X]$ of degree $t$,
        different from $\X$ and $\X-1$.

\medbreak
    \item\label{item:twisted-conjugacy} \emph{Non-trivial (twisted) conjugacy classes in simple groups}.

\medbreak
    \item\label{item:twisted-homogeneous} \emph{Simple twisted homogeneous racks}, that is twisted conjugacy classes of type $(G,u)$, where
    \begin{itemize}
      \item $G = L^t$, with $L$ a simple non-abelian group and $1 < t\in \N$,
      \item $u\in \Aut (L^t)$ acts by
        $$u(\ell_1,\ldots,\ell_t)=(\theta(\ell_t),\ell_1,\ldots,\ell_{t-1}), \quad \ell_1,\ldots,\ell_t \in L,$$
        for some $\theta\in\Aut(L)$.   Furthermore, $L$
        and $t$ are unique, and $\theta$ only depends on
        its conjugacy class in $\Out (L^t)$.
    \end{itemize}
\end{enumerate}

\subsubsection*{Notation.} A \emph{simple rack of type $(L,t,\theta)$} is
a twisted homogeneous rack as in \eqref{item:twisted-homogeneous}.

\subsection{The approach by subracks}\label{subsec:subracks}

\

The experience shows that the following strategy is useful to approach the study of Nichols algebras over finite groups.
However, there are racks that can not be treated in this way.

\begin{strategy}\label{strategy:racks} Let $X$ be a simple rack.

\begin{description}
  \item[Step 1] Is $X$ of type D? In the affirmative, we are done: $X$
and any indecomposable rack $Z$ that admits a rack epimorphism $Z\to X$
collapse, in the sense of Definition \ref{def:rack-collapses}.

\medbreak\item[Step 2]  If not, look for the  abelian subracks of $X$.
  For an abelian subrack $S$ and any $\q\in H^2(X, \ku^{\times})$, look at the
diagonal braiding with matrix $(\q_{ij})_{i,j \in S}$.
  If the Nichols algebra associated to this diagonal braiding has infinite dimension,
and this is known from \cite{H-all}, then so has $\toba(X,\q)$.
  Here you do not need to know all the abelian subracks,
just to find one with the above condition.

\medbreak\item[Step 3] Extend the analysis of Step 2 to cocycles of arbitrary degree.

\medbreak\item[Step 4] Extend the analysis of Steps 2 and 3 to
indecomposable racks $Z$ that admit a rack epimorphism $Z\to X$.
\end{description}
\end{strategy}

The following algorithm is the tool needed to deal with Step 1,
when the rack $X$ is realized as a conjugacy class.

\begin{alg}\label{algo:racks-claseD} Let $\varGamma$ be a
finite group and let $\oc$ be a conjugacy class. Fix $r\in \oc$.

\begin{enumerate}
    \item For any $s\in \oc$, check if $(rs)^2 \neq (sr)^2$; this
    is  equivalent to \eqref{eqn:hypothesis-subrack}.

    \smallbreak
    \item If such $s$ is found, consider the subgroup $H$ generated
    by $r$, $s$. If $\oc^H_r \cap \oc^H_s = \emptyset$, then
    $Y = \oc^H_r \coprod \oc^H_s$ is the decomposable subrack
    we are looking for and $\oc$ is of type D.
\end{enumerate}
\end{alg}

In practice, we implement this algorithm in a recursive way,
running over the maximal subgroups, see \cite{AFGV-espo} for
details.

\bigbreak

Let $X$ be a rack and $S$ a subset of $X$. We denote
$\ll S\gg :=  \bigcap_{\substack{Y \text{subrack},\\ S \subset Y\subset X }}Y.$
If $X$ is a subrack of a group $G$ and $H = \langle S \rangle$, then
$\ll S\gg = \bigcup_{s\in S} \oc^H_s$.

\medbreak
There are racks that could not be dealt with the criterium of type D.

\begin{definition}\label{def:rack-tipoM}
An indecomposable finite rack $X$ is \emph{of type M} if\footnote{M stands for Montevideo, where this notion was discussed by two of the authors.}
for any $r, s\in X$, $\ll\hspace{-2pt} \{r,s\}\hspace{-2pt}\gg$ either is indecomposable or else equals $\{r,s\}$.
\end{definition}

There are racks such that all proper subracks are abelian; for instance, the conjugacy class of type $(2,3)$ in $\sco$
  (here, all proper subracks have at most 2 elements).
More examples of racks  of type M can be found in \cite[Remark 4.2]{AFGV}.

\section{Tools for cocycles}\label{sec:tools}

\subsection{The enveloping group}\label{subsec:enveloping}

\

The enveloping group $\Gx := \langle e_x:
x\in X\vert e_x\, e_y = e_{x\trid y}\, e_x, x,y \in X\rangle$ was introduced in \cite{Bk, fr,
jo}; it was also considered in \cite{lyz, ess, s}. The map $e:X\to \Gx$,
$x\mapsto e_x$ has a universal property:

\medbreak \emph{If $H$ is a group and $f:X\to H$ is a morphism of racks,
then there is a unique morphism of groups $F: \Gx \to H$ such that
$F(e_x) = f_x$, $x\in X$.}

\medbreak In other words, $X\rightsquigarrow \Gx$ is the
adjoint of the forgetful functor from groups to racks.

Since $\varphi: X \to \Sim_{X}$ is a morphism of racks,
there is a group morphism $\Phi: \Gx \to \Sim_X$.
The image, resp. the kernel, of $\Phi$ is denoted $\gax$
(the group of inner automorphisms), resp. $\Gamma_X$ (the defect group).

\medbreak
The group $\gax$ is not difficult to compute in the case of our interest.
As for the defect group,
some properties were established by Soloviev.

\begin{theorem}\label{obs:enveloping-group}
\renewcommand{\theenumi}{\alph{enumi}}   \renewcommand{\labelenumi}{(\theenumi)}
\begin{enumerate}
\item\label{b} If $X$ is a subrack of a group $H$, then $\gax\simeq C/Z(C)$, where $C$
is the subgroup generated by $X$ \cite[Lemma 1.9]{AG-adv}.

\item\label{c} The defect group $\Gamma_X$ is central in
$\Gx$ \cite[Th. 2.6]{s}. Hence $\Gamma_X = Z(\Gx)$ if $\gax$
is centerless.

\item\label{d} The rank of $\Gamma_X$ is the number of
$\gax$-orbits in $X$ \cite[Th. 2.10]{s}. \qed
\end{enumerate}
\end{theorem}

The difficult part of the calculation of the defect group is to
compute its torsion.

\begin{prop}\label{prop:enveloping-group-sm}
Let $s_i = (i \, i+1) \in \sm$. Let $X$ be the rack of
transpositions in $\sm$; this is the conjugacy class of $s_1$.
The enveloping group $\Gx$ of $X$ is a central extension
\begin{equation}\label{eq:enveloping-group-sm}
\xymatrix{0 \ar[0,1] &   \Z \ar[0,1] & \Gx \ar[0,1] & \sm \ar[0,1]
& 0 }
\end{equation}
\end{prop}

\pf By property \eqref{b} above, $\gax\simeq
\sm$. We have to compute $\Gamma_X$. Let $\mathbb B_m$ be the
braid group  and
     as in Subsection \ref{subsec:nichols}; let $\mathbb P_m = \ker \pi$, the pure braid
group. We claim that there is a morphism of groups $\Psi: \mathbb
B_m \to \Gx$ with $\Psi(\sigma_i) = e_{s_i}$, $1\le i \le m-1$. To
prove the claim, we verify the defining relations of the braid
group:
\begin{align*}
\text{If }\vert i-j\vert &\geq 2, \text{ then  }&  e_{s_i}e_{s_j}
&= e_{s_i\trid s_j}e_{s_i} = e_{s_j}e_{s_i} ;
\\
\text{if }\vert i-j\vert & = 1, \text{ then  }&
e_{s_i}e_{s_j}e_{s_i} &= e_{s_i}e_{s_j\trid s_i}e_{s_j} =
e_{s_i\trid(s_j\trid s_i)}e_{s_i}e_{s_j} = e_{s_j}e_{s_i}e_{s_j}
\end{align*}
since $s_i\trid(s_j\trid s_i) = s_j$ in $\sm$. In other words, we
have a commutative diagram
\begin{equation*}
\xymatrix{\mathbb B_m \ar[1,1]^{\pi} \ar[0,2]^-{\Psi} & &  \Gx
\ar@{>>}[1,-1]^{\Phi}
\\ &\sm. &}
\end{equation*}
Clearly, $\Psi$ is surjective and $\ker \Phi = \Psi(\mathbb P_m)$.
Let now $H$ be a group and $f:X\to H$ a morphism of racks. If
$x,y\in X$, then $f_x^2f_y = f_xf_{x\trid y}f_x = f_yf_x^2$ and consequently
$f_{y\trid x}^2 = f_yf_x^2 f_y^{-1} = f_x^2$. Hence for all $x,y
\in X$,
\begin{equation}\label{eq:rels-rack-sm}
f_y^2 = f_x^2 \quad\text{is central in the subgroup generated by }
f(X).
\end{equation}
We call $z = e_{s_i}^2$; this is a central element of $\Gx$ and
does not depend on $i$. Now $\mathbb P_m$ is generated by
$\tau_{ij} = \sigma_j \trid (\sigma_{j+1} \trid (\sigma_{j+2}\trid\dots
\trid(\sigma_{i-1}\trid\sigma_i^{2}))$ for all $j<i$,
see \cite[page 119]{Ar},\cite{Bi}. Hence $\ker \Phi =
\Psi(\mathbb P_m)$ is generated by $\Psi(\tau_{ij}) = z$.

Let now $V$ be a vector space with a basis $(v_x)_{x\in X}$ and
let $q\in \ku$ be a root of 1 of arbitrary order $M$. Define
$f_y\in \GL(V)$ by $f_y(v_x) = q v_{y\trid x}$, $x,y\in X$. Then
$f_xf_y = f_{x\trid y}f_x$ and $f_x^2 = q^2\id$ for any $x,y\in
X$; thus we have a map $F: \Gx\to \GL(V)$ and $F(z)= q^2\id$. This
implies that $z$ is not torsion and the claim is proved.
\epf

\begin{obs}
(i). By \cite[Ch IV, \S 1, no. 1.9, Prop. 5]{B}, there is a section
of sets $T: \sm \to \Gx$ such that $T(ww') = T(w) T(w')$ when
$\ell(ww') = \ell(w) \ell(w')$. Thus the central extension
corresponds to the cocycle $\eta: \sm \times \sm \to \Z$, $\eta(w,w')
= T(w) T(w') T(ww')^{-1}$, $w,w'\in \sm$.

(ii). The proof shows the centrality of $\Gamma_X$ directly
without referring to Theorem \ref{obs:enveloping-group} \eqref{c}.
By Theorem \ref{obs:enveloping-group} \eqref{d}, $z$ is not
torsion; the last paragraph of the proof avoids appealing to this
result.
\end{obs}

\bigbreak
Let $(X,\trid)$ be a rack, $\q:X\times X\to\GL(n,\ku)$
a 2-cocycle of degree $n$ and $(V,c) = (\ku X\otimes\ku^{n}, c^\q)$, \emph{cf.} \eqref{eqn:trenza}.
We discuss how to realize $(V, c)$
as a Yetter-Drinfeld module over a group algebra.
Let $x\in X$ and define $g_x:V\to V$ by
\begin{equation}\label{eqn:gi}
g_x(e_yw)=e_{x\trid y}\q_{x,y}(w), \qquad y\in X, w\in \ku^{n},\end{equation} and
let $\ixq$ be the subgroup of $\GL(V)$ generated by the $g_x$'s, $x\in
X$. Then $g_xg_y = g_{x\trid y}g_x$ for any $x,y\in X$, and $(V, c)$ is a Yetter-Drinfeld module over the group algebra
of $\ixq$, with the natural action and coaction $\delta(e_xv) = g_x\otimes e_xv$, $x\in X$, $v\in \ku^{n}$.

\begin{lema}
Let $F$ be a group provided with:
\begin{itemize}
    \item a group homomorphism $p:F\to \ixq$;
    \item a rack homomorphism $s:X\to F$ such that $p(s_x) = g_x$ and $F$ is generated as a group by
    $s(X)$.
\end{itemize}
Then $(V, c)\in \ydf$, with the action induced by $p$ and coaction
$\delta(e_xv) = s_x\otimes e_xv$, $x\in X$, $v\in \ku^{n}$.  In
particular,  $(V, c)\in \ydgx$.
\end{lema}

\pf If $x,y\in X$ and $w\in \ku^n$, then $\delta(s_x\cdot e_yw) =
\delta(e_{x\trid y}\q_{x,y}(w)) = s_{x\trid y}\otimes e_{x\trid
y}\q_{x,y}(w) = s_{x}s_{y}s_{x}^{-1} \otimes e_{x\trid
y}\q_{x,y}(w)=s_{x}s_{y}s_{x}^{-1} \otimes s_x\cdot w$. Since $F$
is generated by $s(X)$, it follows that $\delta(f\cdot e_yw) =
fs_{y}f^{-1} \otimes f\cdot w$, for all $f\in F$.\epf

As a consequence, the Nichols algebra of the braided vector space
$(\ku X, c_q)$ bears a $\Gx$-grading, that we shall call the
\emph{principal} grading, as opposed to the natural $\N$-grading.
Indeed, if $X$ is abelian, then $\Gx \simeq \Z^{\card X}$ and the principal grading
coincides with the one considered \emph{e.~g.} in \cite{AHS}.

\subsection{The rack cohomology group
$H^2(X, \ku^{\times})$}\label{subsec:H2}

\

We now state some general facts about the cocycles on a rack $X$
with values in the abelian group
$\ku^{\times}$. In this case, the $H^2$ is part of a cohomology theory,
see \cite{AG-adv} and references therein.
An alternative description of $H^2(X, \ku^{\times})$
was found in \cite{EG} through
the enveloping group. Namely,
let $\Fun(X, \ku^{\times})$ be the space of all functions from
$X$ to $\ku^{\times}$ with  right $\Gx$-action
given by $(f\cdot e_x) (y) = f(x \trid y)$, $f: X \to A$, $x,y\in X$.

\begin{lema}\label{lem:EG-H2} \cite{EG}
$H^2(X,\ku^{\times}) \simeq H^1(\Gx, \Fun(X, \ku^{\times}))$. \qed
\end{lema}

In principle, the cohomology of $\Gx$ could be studied via
the Hochschild-Serre sequence from that of $\gax$ and
$\Gamma_X$. However, the computation of the defect group
seems to be very difficult.
There is also a homology theory of racks, related to the computation
we are interested in by the following result.

\begin{lema}\label{lem:AG-47}
\cite[Lemma 4.7]{AG-adv}
$H^2(X,\mathbb{C}^\times)\simeq\textrm{Hom}
(H_2(X,\mathbb{Z}),\mathbb{C}^\times)$. \qed
\end{lema}

\medbreak There is a monomorphism
$\ku^{\times} \hookrightarrow H^2(X,\ku^{\times})$,
since any constant function is a cocycle.
A natural question is to compute the quotient
$H^2(X,\ku^{\times}) / \ku^{\times}$.
\emph{Assume that $X$ is indecomposable}.
If $\q\in Z^{2}(X,\ku^{\times})$, then
\begin{equation}\label{eqn:cociclo-enla-diagonal}
\q_{ii} =  \q_{jj}, \quad \text{for any } i,j\in X.
\end{equation}
Note also that $\q \sim \widetilde{\q}$ as in
\eqref{eqn:rel-cociclos} implies that $\q_{ii} =  \widetilde{\q}_{ii}$ for all $i\in X$.
Therefore the question can be rephrased as follows.

\begin{question}\label{que:cociclo-diagonal-menosuno}
Compute all cocycles $\q\in Z^{2}(X,\ku^{\times})$ such that $\q_{ii} =  -1$.
\end{question}

\subsection{The program \textsf{RiG}}\label{subsec:rig}

\

A program for calculations with racks, that in particular computes the
rack-(co)homology groups, was developed in \cite{GV}. We use it
to compute some cohomology groups of simple racks
that are not of type D, see Theorems \ref{th:racks-an-liquidados}
and \ref{th:racks-sn-liquidados}.

We say that $\sigma\in\s_n$ is of type
$(1^{n_1},2^{n_2},\dots,m^{n_m})$ if the decomposition of $\sigma$
as product of disjoint cycles contains $n_j$ cycles of length $j$,
for every $j$, $1\leq j \leq m$.

\begin{prop}\label{prop:rig} Let $\sigma\in \mathbb S_m$ be of type $(1^{n_1}, 2^{n_2}, \dots,
m^{n_m})$ and let
$$\Oc = \begin{cases}\text{the conjugacy class of
$\sigma$ in $\sm$}, &\text{if } \sigma\notin \am,\\  
\text{the conjugacy class of  $\sigma$ in $\am$}, &\text{if }
\sigma\in \am.
\end{cases}$$
\renewcommand{\theenumi}{\alph{enumi}}   \renewcommand{\labelenumi}{(\theenumi)}
\begin{enumerate}
\item If $m= 5$ and the type is $(2,3)$,
then $H^2(\Oc, \ku^{\times}) = \ku^{\times}\times\G_6$.
\item If $m= 6, 7, 8$ and the type is $(1^n,2)$,
then $H^2(\Oc, \ku^{\times}) = \ku^\times\times\G_2$.
\item If $m= 6$ and the type is $(2^3)$,
then $H^2(\Oc, \ku^{\times}) = \ku^\times\times\G_2$.
\item If $m= 5$ and the type is $(1^2,3)$,
then $H^2(\Oc, \ku^{\times}) = \ku^\times\times\G_6$.
\item If $m= 6$ and the type is $(1,2,3)$,
then $H^2(\Oc, \ku^{\times}) = \ku^\times\times\G_{3}\times\G_6$.

\end{enumerate}
\end{prop}

\begin{table}[ht]
\begin{center}
\caption{Some homology groups of conjugacy classes in $\sm$.}
\label{tab:H_2}
\begin{tabular}{|c|c|c|}
\hline
\multicolumn{2}{|c|}{type of $X$}   & $H_{2}(X,\mathbb{Z})$\tabularnewline
\hline
$\mathbb{S}_{5}$ & $(1\,2)(3\,4\,5)$ & $\mathbb{Z}\oplus\mathbb{Z}/{6}$\tabularnewline
\hline
$\mathbb{A}_{5}$ & $(1\,2\,3)$ & $\mathbb{Z}\oplus\mathbb{Z}/{6}$\tabularnewline
\hline
$\mathbb{S}_{6}$ & $(1\,2)(3\,4)(5\,6)$ & $\mathbb{Z}\oplus\mathbb{Z}/{2}$
\tabularnewline
\hline
$\mathbb{S}_{6}$ & $(1\,2)$ & $\mathbb{Z}\oplus\mathbb{Z}/{2}$\tabularnewline
\hline
$\mathbb{A}_{6}$ & $(1\,2\,3)$ & $\mathbb{Z}\oplus\mathbb{Z}/{3}\oplus\mathbb{Z}/{6}$
\tabularnewline
\hline
$\mathbb{S}_{7}$ & $(1\,2)$ & $\mathbb{Z}\oplus\mathbb{Z}/{2}$\tabularnewline
\hline
$\mathbb{S}_{8}$ & $(1\,2)$ & $\mathbb{Z}\oplus\mathbb{Z}/{2}$\tabularnewline
\hline
\end{tabular}
\par\end{center}
\end{table}

\pf We use \textsf{GAP} and \textsf{RiG} to compute the homology
groups $H_2(\Oc, \Z)$. These results are listed in Table
\ref{tab:H_2}. Now assume  $X$ is a rack and that there exists
$m\in\mathbb{N}_{\geq2}$ such that
$H_2(X,\mathbb{Z})\simeq\mathbb{Z}\oplus\mathbb{Z}/{m_{1}}\oplus
\cdots \oplus \mathbb{Z}/{m_{r}} $. By Lemma \ref{lem:AG-47}, we
have $H^2(X,\mathbb{C}^\times)\simeq
\textrm{Hom}(\mathbb{Z}\oplus\mathbb{Z}/{m_{1}}\oplus \cdots
\oplus \mathbb{Z}/{m_{r}}, \mathbb{C}^\times)\simeq
\mathrm{Hom}(\mathbb{Z},\mathbb{C}^\times)
\times\mathrm{Hom}(\mathbb{Z}/{m_{1}},\mathbb{C}^\times)
\times\cdots \times
\mathrm{Hom}(\mathbb{Z}/{m_{r}},\mathbb{C}^\times)
\simeq\mathbb{C}^\times\times\mathbb{G}_{m_{1}}\times \cdots
\times \mathbb{G}_{m_{r}}$. \epf

\subsection{Twisting}\label{subsec:twisting}

\

There is a method, called twisting,  to deform the
the multiplication of a Hopf algebra \cite{DT};
it is formally dual to the twisting of the comultiplication
\cite{Dr2,Re}.  The relation  with bosonization was established in
\cite{MO}. Here we show how to relate two cocycles over a rack $X$ by a twisting,
in a way that the corresponding Nichols algebras are preserved.

\medbreak
Let $\Hc$ be a Hopf algebra. Let
$\phi:  \Hc\otimes \Hc\to \ku$ be an invertible (with respect to
the convolution) linear map and define a new product by $
x\cdot_{\phi}y = \phi(x_{(1)}, y_{(1)}) x_{(2)} y_{(2)}
\phi^{-1}(x_{(3)}, y_{(3)})$, $x,y\in \Hc$. If $\phi$ is a unitary
2-cocycle, that is for all $x,y,z\in \Hc$,
\begin{align}\label{cociclo-conduno} \phi(x_{(1)}\otimes y_{(1)}) \, \phi(x_{(2)} y_{(2)}\otimes z) &=
\phi(y_{(1)}\otimes z_{(1)}) \, \phi(x\otimes y_{(2)}z_{(2)}),
\\ \label{cociclo-conddos}
\phi(x\otimes 1) &= \phi(1\otimes x) = \varepsilon(x),
\end{align}
then $\Hc_{\phi}$ (the same coalgebra but with multiplication
$\cdot_{\phi}$) is a Hopf algebra.

\begin{theorem}\label{th:mo2} \cite[2.7, 3.4]{MO}
Let $\phi:  \Hc\otimes \Hc\to \ku$ be an invertible unitary
2-cocycle.
\renewcommand{\theenumi}{\alph{enumi}}   \renewcommand{\labelenumi}{(\theenumi)}
\begin{enumerate}
    \item There exists an equivalence of braided categories $\T_{\phi} :
\ydh\to \ydhs$, $V\mapsto V_{\phi}$, which is the identity on
the underlying vector spaces, morphisms and coactions, and
transforms the action of $\Hc$ on  $V$ to $\cdot_{\phi}:
\Hc_{\phi}\otimes V_{\phi}\to V_{\phi}$,
$$
h\cdot_{\phi} v = \phi(h_{(1)}, v_{(-1)}) (h_{(2)}\cdot
v_{(0)})\_0\, \phi^{-1}((h_{(2)}\cdot v_{(0)})\_{-1}, h_{(3)}),
$$
$h\in \Hc_{\phi}$, $v\in V_{\phi}$. The monoidal structure on
$\T_{\phi}$ is given by the natural transformation $b_{V,W}:
(V \otimes W)_{\phi} \to V_{\phi} \otimes W_{\phi}$
\begin{align*}
b_{V,W}(v\otimes w) = \phi(v_{(-1)}, w_{(-1)}) v\_0 \otimes
w\_0, \quad v\in V, w\in W.
\end{align*}
\item $\T_{\phi}$ preserves Nichols algebras: $\toba(V)_{\phi}
\simeq \toba(V_{\phi})$ as objects in $\ydhs$. In particular,
the Poincar\'e series of $\toba(V)$ and $\toba(V_{\phi})$ are
the same. \qed
\end{enumerate}
\end{theorem}

 Let us recall the argument for (b). The functor $\T_\phi$
preserves the braidings; that is, if $c$, resp. $c_\phi$, is the
braiding in $\ydh$, resp. $\ydhs$, then the following diagram
commutes:
\begin{equation}\label{eqn:def-twisted-braiding}
\xymatrix{(V \otimes W)_\phi
\ar@{->}^{\T_\phi(c)}[0,2]\ar@{->}[1,0]_{b_{V,W}}&
& (W \otimes V)_\phi\ar@{->}[1,0]^{b_{W,V}} \\
V_\phi \otimes W_\phi \ar@{->}^{c_\phi}[0,2]& & W_\phi \otimes
V_\phi.}
\end{equation}

Since the ideal of relations of a Nichols algebra is the sum of
the kernels of the various quantum symmetrizers, (b) follows
immediately.

\medbreak Let $G$ be a group. If $\Hc = \ku G$, then a unitary
2-cocycle on $\Hc$ is equivalent to a 2-cocycle $\phi\in Z^2(G,
\ku^{\times})$, that is a map $\phi: G\times G\to \ku^{\times}$
such that
\begin{equation}\label{eqn:cocycle-group}
\phi(g,h) \, \phi(gh, t) = \phi(h, t) \, \phi(g, ht)
\end{equation}

and $\phi(g, e) = \phi(e, g) = 1$ for all $g,h,t\in G$.

\medbreak Let $\sigma, \zeta\in G$, $\oc_\sigma$, $\oc_\zeta$
their conjugacy classes, $(\rho,V)\in \widehat{C_{G}(\sigma)}$,
$(\tau,W)\in \widehat{C_{G}(\zeta)}$. For $\nu\in\oc_\sigma$,
$\xi\in\oc_\zeta$, pick $g_\nu, h_\xi\in G$ such that $g_\nu\trid
\sigma = \nu$, $h_\xi\trid \zeta = \xi$.

\begin{lema}\label{lema:tw-gps}
If $\phi\in Z^2(G, \ku^{\times})$, then the braiding $$c_\phi:
M(\oc_\sigma, \rho)_\phi \otimes M(\oc_\zeta, \tau)_\phi \to
M(\oc_\zeta, \tau)_\phi \otimes M(\oc_\sigma, \rho)_\phi$$ is given
by
\begin{equation}\label{eqn:tw-braiding}
c_\phi(g_\nu v\otimes h_\xi w) = \phi(\nu, \xi) \phi^{-1}(\nu\trid
\xi, \nu) \, \nu \cdot h_\xi w \otimes g_\nu v,
\end{equation}
$v\in V$, $w\in W$.
\end{lema}

\pf By \eqref{eqn:def-twisted-braiding}, since $b_{M(\oc_\sigma,
\rho), M(\oc_\zeta, \tau)}(g_\nu v\otimes h_\xi w) = \phi(\nu,
\xi)\, g_\nu v\otimes h_\xi w.$ \epf

Let now $X$ be a subrack of a conjugacy class $\Oc$  in $G$, $q$ a
2-cocycle on $X$ arising from some Yetter-Drinfeld module $M(\Oc,
\rho)$ with $\dim\rho = 1$ and $\phi\in Z^2(G, \ku^{\times})$.
Define $q^\phi: X\times X\to \ku^\times$ by
\begin{equation}\label{eqn:qphi}
q^\phi_{x y} = \phi(x, y) \phi^{-1}(x\trid y, x) \,q_{xy}, \quad
x,y\in X.
\end{equation}

Then Lemma \ref{lema:tw-gps} and Th. \ref{th:mo2} imply that
\begin{equation}\label{eqn:nichols-tw-racks}
\text{The Poincar\'e series of $\toba(X, q)$ and $\toba(X,
q^\phi)$ are equal.}
\end{equation}

\begin{obs}\label{obs:twist-racks} If $X$ is any rack, $q$ a
2-cocycle on $X$ and $\phi: X\times X\to \ku^\times$, then define
$q^\phi$ by \eqref{eqn:qphi}. It can be shown that $q^\phi$ is a
2-cocycle iff
\begin{multline}\label{eqn:condition-tw-racks}
\phi(x,z)\phi(x\trid y, x\trid z)\phi(x\trid(y\trid z),x)\phi(
y\trid z, y)
\\= \phi(y,z)\phi(x, y\trid z)\phi(x\trid(y\trid z), x\trid y)\phi(x\trid z, x)
\end{multline}
for any $x,y,z\in X$. Thus, if $X$ is a subrack of a group $G$ and
$\phi\in Z^2(G, \ku^{\times})$, then $\phi\vert_{X\times X}$
satisfies \eqref{eqn:condition-tw-racks}.

\end{obs}

\begin{definition}\label{defi:twist-racks} The 2-cocycles $q$ and $q'$ on $X$
are \emph{equivalent by twist} if there exists $\phi: X\times X\to \ku^\times$ such that
$q' = q^\phi$ as in  \eqref{eqn:qphi}.
\end{definition}

\section{Simple affine racks}\label{sec:simple-affine}

Let $p$ be a prime, $t\in \N$  and $f\in \fp[\X]$ of degree $t$,
monic irreducible and different from $\X$ and $\X-1$. Let $T$ be
the companion matrix of $f$ and $\Q_{\fp^t, f} := \Q_{\fp^t, T}$
be the associated affine rack; this will be simply denoted by $\Q$
if no emphasis is needed. Alternatively, set $q=p^{t}$ and
identify $\mathbb{F}_q$ with $\mathbb{F}_p^t$. Then the action of
$T$ corresponds to multiplication by $a$, which is the class of
$\X$ in $\fp[\X]/(f)$. Note that $a$ generates $\mathbb{F}_q$ over
$\mathbb{F}_p$.


\begin{question}\label{que:simple-affine-subracks} Find the proper subracks of $\Q$.
\end{question}

We expect that the simple affine racks will have very few subracks.
In fact, they have no abelian subracks with more than one element
\cite[Remark 3.13]{AFGV}.


\begin{prop}
If $a$ generates $\mathbb{F}_q^\times$, then any proper subrack of
$\Q_{\kc,a}$ is trivial.
\end{prop}

\begin{proof}
Let $X$ be a subrack of $\Q_{\kc,a}$ with more than one element.
Let $x,y\in X$ with $x\ne y$. By definition we have
$\varphi_x^n(y)\in X$ for all $n\in\mathbb{N}$. Since
$\varphi_x^n(y)=(1-a^n)x+a^ny$, for all $n\in\mathbb{N}$,  we have
that
\[
A=\{(1-a^n)x+a^ny\mid 0\leq n\leq q-1\}\subseteq X.
\]
Moreover, $A$ has $q$ elements. Indeed, suppose there exist $m\neq
n$ such that $(1-a^n)x+a^ny=(1-a^m)x+a^my$. Then
$x(a^m-a^n)=y(a^m-a^n)$ which implies that $ x=y $, a
contradiction. Since $ A\subseteq  X \subseteq\Q_{\kc,a} $ and $
|\Q_{\kc,a}| =q$ we have that $X=\Q_{\kc,a}$.
\end{proof}

In the particular case $t=1$, we can say more: any proper subrack of an affine rack
with $p$ elements is trivial.

\begin{prop}
Let $1\ne a\in\mathbb{F}_p^\times$. Then any proper subrack of the affine rack
$\Q_{\mathbb{F}_p,a}$  is trivial.
\end{prop}

\begin{proof}
Let $x\ne y$ be two elements of $\mathbb{F}_p$. It is enough to
show that the subrack generated by $x$ and $y$ is $\mathbb{F}_p$.
Let
\[
F_{a,m}(n_1,n_2,...,n_m)=\sum_{j=1}^{m}(-1)^{j+1}a^{n_j+\cdots+n_m}.
\]
Note that $a+aF_{a,2k}(n_1,n_2,...,n_{2k})$ $=
F_{a,2k+1}(n_1,n_2,...,n_{2k},1)$. Then
\begin{align}
\label{eq:par}\varphi_y^{n_{2k}}\varphi_x^{n_{2k-1}}
\cdots\varphi_x^{n_1}(y)&=y+(y-x)F_{a,2k}(n_1,n_2,...,n_{2k}),\\
\label{eq:impar}\varphi_y^{n_{2k+1}}\varphi_x^{n_{2k}}
\cdots\varphi_x^{n_1}(y)&=x+(y-x)F_{a,2k+1}(n_1,n_2,...,n_{2k+1}).
\end{align}
Let $z\in\mathbb{F}_p$, then
\begin{gather}
\label{eq:Fp}
z=\varphi_y^{n_{2k}}\varphi_x^{n_{2k-1}}\cdots\varphi_x^{n_1}(y)
\end{gather}
has at least one solution. In fact, let $n_j=(-1)^j$. Equation
\eqref{eq:par} implies that \eqref{eq:Fp} can be re-written as
$z=y+(y-x)(1-a)k$. Then the result follows by taking
$k=(z-y)(1-a)^{-1}(y-x)^{-1}$.
\end{proof}


\section{Conjugacy classes in non-abelian simple groups}\label{sec:type1-theta-id}

\subsection{Alternating groups}\label{subsect:alt}

\


\begin{theorem}\label{th:racks-an-liquidados} \cite[Th. 4.1]{AFGV}
Let $\sigma\in \am$, $m\geq 5$. If the type of $\sigma$ is NOT any
of $(3^2)$; $(2^2,3)$; $(1^n,3)$; $(2^4)$; $(1^2,2^2)$; $(1,2^2)$;
$(1,p)$, $(p)$ with $p$ prime, then the conjugacy class of
$\sigma$ in $\am$ is of type D.\qed
\end{theorem}

\subsection{Sporadic groups}\label{subsect:sporadic}

\begin{theorem}\label{th:racks-liquidados} \cite{AFGV-espo, logbook}
  If $G$ is a sporadic simple group and  $\oc$ is a non-trivial conjugacy class of
  $G$ NOT listed in Table \ref{tab:0}, then $\oc$ is of type D.
  \qed
\end{theorem}

\begin{table}[ht]
\begin{center}
\caption{Conjugacy classes of sporadic groups which are not known to be of type D;
those which are NOT of type D appear in bold.}\label{tab:0}
\begin{tabular}{|p{1cm}|c||p{1cm}|c|}
\hline $G$ & {\bf Classes} & $G$ & {\bf Classes}
\\ \hline  $M_{11}$ &  \textup{\bf 8A, 8B, 11A, 11B} &
$Co_{1}$ &  \textup{3A, 23A, 23B}
\\ \hline  $M_{12}$ &  \textup{\bf 11A, 11B} &
$J_{1}$ &  \textup{\bf 15A, 15B, 19A, 19B, 19C}
\\ \hline  $M_{22}$ &  \textup{\bf 11A, 11B} &  $O'N$ &  \textup{31A, 31B}
\\ \hline  $M_{23}$ &  \textup{\bf 23A, 23B} &
$J_{3}$ &  \textup{5A, 5B, 19A, 19B}
\\ \hline  $M_{24}$ &  \textup{\bf 23A, 23B} &    $Ru$ &  \textup{29A, 29B}
\\ \hline  $J_{2}$ &  \textup{\bf 2A, 3A} &    $He$ &  \textup{all of type D}
\\ \hline  $Suz$ &  \textup{3A} &  $Fi_{22}$ &  \textup{{\bf 2A}, 22A,
22B}
\\ \hline  $HS$ &  \textup{11A, 11B} &   $Fi_{23}$ &  \textup{{\bf 2A},  23A, 23B}
\\ \hline  $McL$ &  \textup{11A, 11B} &  $HN$ &  \textup{all of type D}
\\ \hline  $Co_{3}$ &  \textup{23A, 23B} & $Th$ & \textup{all of type D}
\\ \hline  $Co_{2}$ &  \textup{{\bf 2A}, 23A, 23B} &  $T$ &  \textup{2A}
\\ \hline
\end{tabular}

\vspace{5pt}

\begin{tabular}{|p{0,8cm}|c|}
\hline $G$ & {\bf Classes}
\\ \hline  $Ly$ &  \textup{33A, 33B, 37A,
37B, 67A, 67B, 67C}
\\ \hline   $J_4$ &  \textup{29A, 37A, 37B, 37C, 43A, 43B, 43C}
\\ \hline   $Fi'_{24}$ &  \textup{23A, 23B, 27B, 27C, 29A, 29B, 33A, 33B, 39C, 39D}
\\ \hline $B$ &  \textup{2A, 16C, 16D, 32A, 32B, 32C, 32D,
34A, } \\ \hline & \textup{46A, 46B, 47A, 47B}
\\ \hline $M$ &  \textup{32A, 32B, 41A,
46A, 46B, 47A, 47B, 59A, 59B, } \\ \hline & \textup{69A,
69B, 71A, 71B, 87A, 87B, 92A, 92B, 94A, 94B}
\\ \hline\end{tabular}
\end{center}
\end{table}

\subsection{Finite groups of Lie type}\label{subsect:-finite-lie}

\

Let $p$ be a prime number, $m\in \N$ and $q =p^m$. Let $\G$ be an
algebraic reductive group defined over the algebraic closure of
$\kc$ and $G = \G(\kc)$ be the finite group of $\kc$-points. Let
$x\in G$; we want to investigate the orbit $\Oc_x^G$ of $x$ in
$G$. Let $x = x_sx_u$ be the Chevalley-Jordan decomposition of $x$
in $\G$; then $x_s, x_u\in G$. Let $\Kb = C_{\G}(x_s)$, a
reductive subgroup of $\G$ by \cite[Thm. 2.2]{Hu},
and let $\Lb$ be its semisimple part;
then $K := \Kb\cap G = C_{G}(x_s)$, by \cite[Prop. 9.1]{Bo}.
Since $x_u\in K$, we conclude from Subsection
\ref{subsec:racks-typeD}  that
$$
\Oc_{x_u}^K \text{ is a subrack of } \Oc_x^G.
$$
Therefore, we are reduced to investigate the orbits $\Oc_x$
when $x$ is either semisimple (the case $x = x_s$)
or unipotent (by the reduction described).

The first step of the Strategy proposed in
Subsection \ref{subsec:subracks} consists of finding
subracks of type D of conjugacy classes of
semisimple or unipotent elements.
We believe that most semisimple conjugacy classes are of type D.
We give now some evidence for this belief,
using techniques with involutions and elements of
a Weyl group associated to a fixed $\kc$-split torus.

\medbreak
Let $n>1$, $\xi \in \kc^{\times}$ so that $\ord \xi = m $ divides
$q-1$ and $a \in \kc^{\times}$. For all  $x = (x_{1},\ldots, x_{n}) \in
(\Z/m)^{n}$ such that $\sum_{i=0}^{n} x_{i} \equiv 0 \pmod
m$ define $n_{a}$ to be the companion matrix of the polynomial
$X^{n}-a$, $\xi_{x}= \diag(\xi^{x_{1}},\ldots, \xi^{x_{n}})$  and $\mu_{x}= n_{a}\xi_{x}$.

$$
\mu_{x} = \left( \begin{matrix}
         0 & & \ldots  & 0 & a\xi^{x_{n}}\\
         \xi^{x_{1}} & 0 & \ldots & 0 & 0\\
         0 & \xi^{x_{2}} & \ldots & 0 & \vdots \\
       \vdots & & \ddots & \vdots & 0\\
         0 &\ldots & \ldots & \xi^{x_{n-1}} & 0
              \end{matrix}\right)  \in \GL(n,\kc).
$$

Let $X_{a,\xi}= \{\mu_{x}:\ \sum_{i=1}^{n}
x_{i} \equiv 0\pmod m\}$, a subset of the conjugacy class of
$n_{a}$ in $\mathbf{GL}(n,\kc)$
(that is, the set of
matrices with minimal polynomial $T^n -a$). If $a= -1$, then
$X_{a,\xi} \subseteq \SL(n,\kc)$.

The following proposition is a generalization of
\cite[Example 3.15]{AF3}.

\begin{prop}\label{prop:rack-afin-lie}
Assume that $(n, q-1)\neq 1$; that
$q>3$,  if $n = 4$; and that $q>5$,  if $n = 2$.
Then the conjugacy class of $n_a$ is of type D.
\end{prop}

\pf
Pick a generator $\xi$ of $\kc^{\times}$. We claim that $X_{a,\xi}$ is a subrack of the conjugacy class of
$n_{a}$ in $\mathbf{GL}(n,\kc)$, isomorphic to the affine rack
$\Q_{(\Z/(q-1))^{n-1},g}$, with $g(x_{1},\ldots ,x_{n-1})= \left(-\sum_{i=1}^{n-1} x_{i},
x_{1},\ldots, x_{n-2}\right)$.
A direct computation shows that $\mu_{x}\mu_{y}\mu_{x}^{-1} =
\mu_{x\trid y}$, with
\begin{align*}
x\trid y & = (x_{1} + y_{n} -x_{n},
x_{2} + y_{1} -x_{1},\ldots, x_{n} + y_{n-1} -x_{n-1}).
\end{align*}
Thus, the map $\varphi: X_{a,\xi}\to \Q_{(\Z/(q-1))^{n-1},g}$ given
by $\varphi(\mu_{x})= (x_{1},\ldots ,x_{n-1})$ is a rack isomorphism and the claim is proved.
The proposition follows now from \cite[Lemma 2.2]{AFGV-thr}, for $n >2$, or \cite[Lemma 2.1]{AFGV-thr}, for $n=2$.
\epf

The conjugacy class of involutions in $\PSL(2,\kc)$ for $q\in\{5,7,9\}$ is not of
type D. For $q>9$
we have the following result.

\begin{cor}
\begin{enumerate}
\item[(a)] Assume that $q\equiv 1\pmod 4$ and $q>9$.  Then the conjugacy class of involutions of
$\PSL(2,\kc)$ is of type D.
\item[(b)] Assume that $q\equiv 3\pmod 4$ and $q>9$.
Then the conjugacy class of involutions of
$\PGL(2,\kc)$ is of type D.
\end{enumerate}
\end{cor}

\pf 
(a) Let $\xi \in \kc$ such that $\kc^{\times}= \langle \xi \rangle$.
By Proposition \ref{prop:rack-afin-lie} with $a = -1$, the
subrack $X= \{\left(\begin{smallmatrix}
     0 & -\xi^{-x}\\ \xi^{x} & 0\end{smallmatrix}\right):
     x\in \Z/(q-1)\}$ of the conjugacy class
of $n_{-1}$ in $\GL(2,\kc)$ is isomorphic to the dihedral rack $\D_{q-1}$.
Let $\pi: \GL(2,\kc) \to \PGL(2,\kc)$ be the canonical projection.
Then $\pi\left(\begin{smallmatrix}
     0 & -\xi^{-x}\\ \xi^{x} & 0\end{smallmatrix}\right) \in \PSL(2,\kc)$ for all
$x\in \Z/(q-1)$ and whence $\pi(X)$
is a subrack of the unique conjugacy class
of involutions in $\PSL(2,\kc)$. Now $\pi\left(\begin{smallmatrix}
     0 & -\xi^{-x}\\ \xi^{x} & 0\end{smallmatrix}\right) = \pi\left(\begin{smallmatrix}
     0 & -\xi^{-y}\\ \xi^{y} & 0\end{smallmatrix}\right)$ iff $\xi^{x} = - \xi^{y}$,
     hence $\pi(X) \simeq\D_{(q-1)/2}$, which is of type D if $(q-1)/2$ is even and $>4$.

\medbreak
(b) Let $L =\{\left(\begin{smallmatrix} a & b\\ -b & a\end{smallmatrix}
\right):\ a,b\in \kc \}$. Then $L$ is a quadratic field extension of $\kc$
and $|L|=q^{2}$. Consider now the group map $\det: L^{\times}
\to \kc^{\times}$ given by the determinant. 
Since every element in
a finite field is a sum of squares, the 
kernel is a subgroup
of $L^{\times}$ of order 
 $\frac{|L^{\times}|}{|\kc^{\times}|} = \frac{q^{2}-1}{q-1}$. Since
$L^{\times}$ is cyclic, there exist $a, b\in \kc$ such that $a^{2} +
b^{2} = 1$ and $ \theta= \left(\begin{smallmatrix} a & b\\ -b & a\end{smallmatrix}
\right)$ generates $\ker \det$, \textit{i.e.} it has order $q+1$.  
Note that, as $q\equiv 3 \pmod 4$, 
$\theta$ is contained in a non-split torus.  

Let $n = \left(\begin{smallmatrix} 0 & 1\\ 1 & 0\end{smallmatrix}
\right)$. Then  
the
subrack $X= \{\mu_{x}=
     n\theta^{x}:
     x\in \Z/(q+1)\}$ of the conjugacy class
of $n$ in $\GL(2,\kc)$ is isomorphic to the dihedral rack $\D_{q+1}$.
Taking $\pi$ as in (a), we have that 
$\pi(X) \simeq\D_{(q+1)/2}$
is a subrack of the unique conjugacy class
of involutions in $\PGL(2,\kc)$, which is of type D if $(q+1)/2$ is even and $>4$.
\epf

A similar argument as in the proof of proposition \ref{prop:rack-afin-lie}
applies with weaker hypothesis to matrices whose rational form contains $n_a$.

\begin{prop} Let $x \in \mathbf{GL}(N,\kc)$ be a semisimple element
 whose rational form $x$
is $\left(\begin{smallmatrix}
     n_{a} & 0\\
0& B_{1}
    \end{smallmatrix}\right)$. Suppose
there exists $B_{2} \in \mathbf{GL}(N-n,\kc)$ such that
$B_{2}\neq B_{1}$, $B_{2}\sim
B_{1}$ and $B_{1}B_{2}=B_{2}B_{1}$.
Then the conjugacy class of $x$ is of type D for all $n \neq 2, 4$;
or $n=4$ and $q>3$; or $n=2$ and $q$ is odd.
\end{prop}

\pf Pick a generator $\xi$ of $\kc^{\times}$ and let $\mu_x$ be as above.
Let $$X_{i}=\big\{\big(\begin{smallmatrix}
     \mu_{x} & 0\\
     0& B_{i}    \end{smallmatrix}\big):\ \sum_{j=1}^{n}x_{j}\equiv 0 \pmod {q-1}\big\},$$
$i=1,2$ and $X=X_{1}\coprod X_{2}$.
Since $\left(\begin{smallmatrix}
     \mu_{x} & 0\\
0& B_{i}
    \end{smallmatrix}\right)\trid \left(\begin{smallmatrix}
     \mu_{y} & 0\\
0& B_{j}
    \end{smallmatrix}\right) =
    \left(\begin{smallmatrix}
     \mu_{x\trid y} & 0\\
0& B_{j}
    \end{smallmatrix}\right)$
and
$X_{1}\cap X_{2} = \emptyset$, we see that
$X$ is a decomposable rack and each $X_i$
is isomorphic to an affine rack,
by the proof of Proposition \ref{prop:rack-afin-lie}.
If $x= (0,\ldots, 0)$, $y=(1,0,\ldots, 0)$,
$s= \left(\begin{smallmatrix}
     \mu_{x} & 0\\
0& B_{1}
    \end{smallmatrix}\right)$ and
    $r=\left(\begin{smallmatrix}
     \mu_{y} & 0\\
0& B_{2}
    \end{smallmatrix}\right)$, then $r\trid(s\trid(r\trid s)) \neq s$, 
    by a straightforward computation, see the proof of \cite[Lemma 2.2]{AFGV-thr},
    whenever the prescribed restrictions on $n$ hold.
\epf

Assume now that $\G$ be a Chevalley group and
denote by $G=\G(\kc)$ the group of $\kc$-points. Let
$T $ be a $\kc$-split torus in $G$ and $W= N_{G}(T)/C_{G}(T)$
the corresponding Weyl group. 
Let $\sigma \in W$ and 
$n_{\sigma} \in N_{G}(T)$ be a representative of $\sigma$.
Since $W$ stabilizes $T$, the adjoint action of 
$n_{\sigma}$ on $T$
defines an automorphism $ g_{\sigma} $ of $ (\Z/(q-1))^{n} $. Indeed,
without loss of generality, we may assume that $T=
\kc^{\times}\times \cdots \times \kc^{\times}$ and
$\kc^{\times} =\langle \xi \rangle$, with $\xi \in \kc^{\times}$.
Then for all $t \in T$, there exists $x \in (\Z/(q-1))^{n} $
such that $t = \xi_{x}=\diag(\xi^{x_{1}},\ldots, \xi^{x_{n}})$, 
$n=\dim T$, and the automorphism
is defined by $n_{\sigma}\xi_{x}n_{\sigma}^{-1}
= \xi_{g_{\sigma}(x)}$. 

The following proposition is a generalization of
Proposition \ref{prop:rack-afin-lie}.

\begin{prop}  Let $\sigma \in W$ and
$n_{\sigma} \in N_{G}(T)$ be a representative of $\sigma$.
Assume there exists $x \in (\Z/(q-1))^{n}$ 
such that $x \notin \Imm	 (\id - g_{\sigma})$ and
$ x-g_{\sigma}(x) +g_{\sigma}^{2}(x) - g_{\sigma}^{3}(x) \neq 0 $. 
Then 
the conjugacy class of $n_{\sigma}$ in $G$ is
of type $D$.
\end{prop}

\pf 
Consider the set
$ X_{\sigma, \xi} = \{\mu_{y}=n_{\sigma}\xi_{y}:\ y\in (\Z/(q-1))^{n} \}$.
Then $X_{\sigma, \xi}$ is a (non-empty) rack isomorphic to
the affine rack $( (\Z/(q-1))^{n} , g_{\sigma})$. Indeed, since
\begin{align*}
\mu_{x}\mu_{y}\mu_{x}^{-1} & =
n_{\sigma}\xi_{x}n_{\sigma}\xi_{y}\xi_{x}^{-1}n_{\sigma}^{-1}=
n_{\sigma}\xi_{x}n_{\sigma}\xi_{y-x}n_{\sigma}^{-1}
=
n_{\sigma}\xi_{x}\xi_{g_{\sigma}(y-x)}\\
&    =
n_{\sigma}\xi_{g_{\sigma}(y)+ (1-g_{\sigma})(x)}
= \mu_{x\trid y},
\end{align*}
the map $\varphi: X_{\sigma, \xi} \to ( (\Z/(q-1))^{n} , g_{\sigma})$
given by $\varphi(\mu_{x}) = x$ defines a rack isomorphism.
Since $x \notin \Imm	 (\id - g_{\sigma})$, 
$ X_{\sigma, \xi} $ contains at least two cosets with respect
to $ \Imm (1-g_{\sigma}) $. If we take	
$s=\mu_{0}$ and $r=\mu_{x}$, then
$$r\trid(s\trid(r\trid s))= \mu_{x\trid(0\trid (x\trid 0))}=
\mu_{x-g_{\sigma}(x) +g_{\sigma}^{2}(x) - g_{\sigma}^{3}(x)},$$
which implies that $X_{	\sigma, \xi}$ is of type $D$.
\epf



\section{Twisted conjugacy classes in simple non-abelian groups}\label{sec:type1-theta-not-id}

In this section we consider twisted conjugacy classes in simple
non-abelian groups defined by non-trivial outer automorphisms.
These can be realized as conjugacy classes in the following way.
Pick a representative of $\theta$ in $\Aut(L)$, called also $\theta$,
and form the semidirect product $L\rtimes \langle\theta \rangle$.
Then the racks of type $(L, 1, \theta)$ are the conjugacy classes
of $L\rtimes \langle\theta \rangle$ contained in $L\times \{\theta
\}$.

\subsection{Alternating groups}\label{subsect:alt-theta-not-id}

\

Since $\am\rtimes \Z/2 \simeq \s_m$, the racks of this type are
the conjugacy classes in $\s_n$ not intersecting $\an$. We keep
the notation from subsection \ref{subsect:alt}.
Assume that $m\geq 5$.

\begin{theorem}\label{th:racks-sn-liquidados} \cite[Th. 4.1]{AFGV}
Let $\sigma\in \sm - \am$. If the type of $\sigma$ is
neither $(2,3)$, nor  $(2^3)$, nor $(1^n,2)$, then the conjugacy
class of $\sigma$ is of type D.\qed
\end{theorem}

Notice that the racks of type $(2^3)$ and $(1^4,2)$ are isomorphic.
As we see, the only example, except for the type $(2,3)$, is $(1^n,2)$.
We treat it in the following Subsection.

\subsection{The Fomin-Kirillov algebras}\label{subsec:fk}

\

Let $X = \Oc^m_2$ be the rack of transpositions in $\sm$, $m \geq 3$.
As shown in \cite{MS}, see also \cite{afz}, there are two cocycles $\q\in Z^{2}(X,\ku^{\times})$ arising from Yetter-Drinfeld modules over $\ku \sm$ and
such that $\q_{ii} =  -1$ for all (some) $i\in X$. These are either $\q = -1$ or
$\q = \chi$, the cocycle given by $\chi(\sigma, \tau) =
\begin{cases}
  1,  & \mbox{if} \ \sigma(i)<\sigma(j) \\
  -1, & \mbox{if} \ \sigma(i)>\sigma(j).
\end{cases}$, if $\tau, \sigma$ are transpositions, $\tau=(ij)$ and
$i<j$. Furthermore, their classes in $Z^{2}(X,\ku^{\times})$ are different.
Hence, we have a monomorphism $\ku^\times\times\G_2 \hookrightarrow H^2(\Oc^m_2, \ku^{\times})$.

\begin{question}\label{que:H2-trasp}
Is $H^2(\Oc^m_2, \ku^{\times}) \simeq \ku^\times\times\G_2$ for $m \geq 4$?
\end{question}

We conjecture that the answer is yes; Proposition \ref{prop:rig}
(b) gives some computational support to this conjecture, and
Proposition \ref{prop:enveloping-group-sm} should be useful for
this.

\medbreak We turn now to the Nichols algebras associated to $X = \Oc^m_2$.

\begin{itemize}
\renewcommand{\labelitemi}{$\diamond$}

\medbreak   \item If $\q\in Z^{2}(X,\ku^{\times})$ arises from  a Yetter-Drinfeld module over $\ku \sm$ and
$\q_{ii} \neq  -1$, then $\dim \toba(X, \q) = \infty$ \cite[Theorem 1]{afz}. In fact, assume that $m\ge 4$. Then it can be shown that
$\dim \toba(X, \q) = \infty$ for any $\q \in \ku^\times\times\G_2 \hookrightarrow H^2(\Oc^m_2, \ku^{\times})$
such that $\q_{ii} \neq  -1$, just looking at the abelian subrack $\{(12), (34)\}$.

\medbreak   \item The Nichols algebras $\toba(\Oc^m_2, -1)$ and $\toba(\Oc^m_2, \chi)$ are finite-dimensional for $m =3,4,5$, see Table \ref{tab:app}.
Indeed, the Hilbert series of $\toba(\Oc^m_2, -1)$ and $\toba(\Oc^m_2, \chi)$ are equal.

\medbreak   \item The quadratic Nichols algebra of a braided vector space $V$ is
$\widehat{\toba}_2(V) = T(V)/ \langle\ker Q_2\rangle $, cf. \eqref{eqn:quantum-symmetrizer}; clearly, here is an epimorphism $\widehat{\toba}_2(V) \to \toba(V)$. The Nichols algebras $\toba(\Oc^m_2, -1)$ and $\toba(\Oc^m_2, \chi)$ are quadratic for $m =3,4,5$. Furthermore, $\toba(\Oc^m_2, \chi)$ appears in \cite{FK} in relation with the quantum cohomology of the flag variety.

\medbreak   \item It is known if the Nichols algebras $\toba(\Oc^m_2, -1)$ and $\toba(\Oc^m_2, \chi)$ are finite-dimensional, nor if
they are quadratic, for $m \geq 6$.
\end{itemize}

\begin{question}\label{que:trasp-sm}
Are the cocycles $-1$ and $\chi$  equivalent by twist?
Recall that $H^{2}(\sm,\ku^{\times}) \simeq \Z/2$ for $m\geq4$, see \cite{schur}.
\end{question}

A positive answer to Question \ref{que:trasp-sm} would explain the similarities between
the Nichols algebras $\toba(\Oc^m_2, -1)$ and $\toba(\Oc^m_2, \chi)$.

\subsection{Sporadic groups}\label{subsect:sporadic-theta-not-id}

\

The sporadic groups with non-trivial outer automorphisms group are
$M_{12}$, $M_{22}$, $J_2$, $Suz$, $HS$, $McL$, $He$, $Fi_{22}$,
$Fi'_{24}$, $O'N$, $J_3$, $T$ and $HN$. For any group $L$ among
these, the outer automorphisms group is $\mathbb Z/2$ and
$\Aut(L)\simeq L\rtimes \Z/2$. Hence we need to consider the
conjugacy classes in $\Aut(L)-L$.

\begin{theorem}\label{th:esporadicos-twisted-liquidados} \cite{FV}
Let $G$ be one of the following sporadic simple groups: $M_{12}$, $M_{22}$, $J_2$,
$Suz$, $HS$, $McL$, $He$, $O'N$, $J_3$ or $T$. If  $\oc$ is the conjugacy class of a non-trivial
element in $\Aut(G)-G$ NOT listed in Table \ref{tab:twisted-spor}, then $\oc$ is of type D.
  \qed
\end{theorem}

\begin{table}[ht]
\begin{center}
\caption{Twisted conjugacy classes which are not known to be of type D}\label{tab:twisted-spor}
\begin{tabular}{|c|c|c|c|c|c|}
 \hline
\textbf{Group} & $\mathrm{Aut}(M_{22})$ & $\mathrm{Aut}(J_{3})$ & $\mathrm{Aut}(HS)$ & $\mathrm{Aut}(McL)$ & $\mathrm{Aut}(ON)$ 
\tabularnewline \hline

\textbf{Classes}  &  2A & 34A, 34B & 2C &  22A, 22B &  38A, 38B, 38C
\tabularnewline \hline
\end{tabular}

\end{center}
\end{table}

The groups $\Aut(Fi_{22})$, $\Aut(Fi_{24}')$ and $\Aut(HN)$ are
being object of present study, see \cite{FV}.



\section{On twisted homogeneous racks}\label{subsec:s-noid}

In this section, we fix a simple non-abelian group $L$, an integer $t > 1$ and
$\theta\in\Out(L)$; by abuse of notation, we call also by $\theta$ a
representative in $\Aut (L)$. The representative of the trivial element is
chosen as the trivial automorphism.
Let $u\in \Aut (L^t)$ act by
        $$u(\ell_1,\ldots,\ell_t)=(\theta(\ell_t),\ell_1,\ldots,\ell_{t-1}), \quad \ell_1,\ldots,\ell_t \in L.$$

The twisted conjugacy class of $(x_{1},\ldots, x_{t}) \in L^{t}$
is called a \emph{twisted homogeneous rack} (THR for short) of class $(L,t,
\theta)$ and denoted  $\C_{(x_{1},\ldots, x_{t})}$. Let also $
\C_{\ell} := \C_{(e,\ldots , e,\ell)}$, $\ell\in L$. The set of
twisted homogeneous racks of class $(L,t, \theta)$ is
parameterized by the set of twisted conjugacy classes of $L$ under
$\theta$ \cite[Prop. 3.3]{AFGV-thr}. Namely,
\begin{enumerate}

\item If  $(x_{1},\ldots, x_{t}) \in L^{t}$ and $\ell = x_{t}x_{t-1}\cdots x_{2}x_{1}$,
then $\C_{(x_{1},\ldots ,x_{t})} = \C_{\ell}$.

\item  $ \C_{\ell}= \C_k$ iff $k\in \oc^{L, \theta}_{\ell}$; hence $$\C_{\ell}=\{(x_{1},\ldots, x_{t})\in L^{t}:\
x_{t}x_{t-1}\cdots x_{2}x_{1}\in \oc^{L, \theta}_{\ell}\}.$$

\end{enumerate}

In \cite{AFGV-thr}, we have developed some techniques
to check whether $\C_{\ell}$ is of type D;
so far, these techniques are more useful in the case $\theta = \id$.
For illustration, we quote:

\begin{itemize}
  \item If $\ell \in L$ is quasi-real of type $j$,
$t\geq 3$ or $t=2$ and $\ord(\ell) \nmid 2(1-j)$, then $\cl$ is of type D.

\medbreak\item If $\ell$ is an involution and $t > 4$ is even, then $\cl$ is of type D.

\medbreak\item If $\ell$ is an involution, $t$ is odd and
$\oc^{L}_\ell$ is of type D, then so is $\cl$.

\medbreak\item If $(t, \vert L\vert)$ is divisible by an odd prime $p$, or if
$(t, \vert L\vert)$ is divisible by $p=2$ and $t \geq 6$, then $\C_e$ is of type D.

\medbreak\item If $L =\aco$ or $\as$  and $t=2$,
then $\C_e$ is not of type D (checked with  \textsf{GAP}).

\end{itemize}

In other words, at least when $\theta = \id$,
the worse cases are either when $\ell$ is an involution and $t =2, 4$,
or else when $\ell = e$.

As an application of these techniques, we have the following result.

\begin{theorem}\label{th:thr} \cite{AFGV-thr}
    Let $L$ be $\an$, $n\ge 5$, or a sporadic group, $t \ge 2$ and $\ell \in L$. If  $\cl$ is a twisted
homogeneous rack of class $(L, t, \id)$ not listed in Tables
\ref{tab:MainThm:An:id} and \ref{tab:MainThm:espo:id}, then $\cl$ is of type D. \qed
\end{theorem}

\begin{table}[th]
\begin{center}
\caption{THR $\cl$ of type $(\A_n,t,\theta)$, $\theta=\id$,
$t\geq2$, $n\geq 5$, which are not known to be of type D. Those not of type D are
in bold.}\label{tab:MainThm:An:id}
\begin{tabular}{|c|c|c|c|} 
\hline  $n$ & $\ell$ & Type of $\ell$ & $t$ \\

\hline   any & $e$ &$(1^n)$ & odd, $(t,n!)=1$\\
   5 & & ${\bf (1^5)}$ & {\bf 2}\\
   5 & & $(1^5)$ & 4\\

   6 & & ${\bf (1^6)}$ & {\bf 2}\\

\hline
5 & involution & $(1,2^{2})$ & 4, odd\\

6 & & $(1^{2},2^{2})$ & odd\\

8 & &$(2^{4})$ & odd\\

\hline
any & order 4& $(1^{r_1},2^{r_2},4^{r_4})$, $r_4>0$, $r_2+r_4$ even & 2\\

\hline
\end{tabular}
\end{center}
\end{table}

\begin{table}[th]
\begin{center}
\caption{THR $\cl$ of type $(L,t,\theta)$, with $L$ a sporadic
group, $\theta=\id$, which are not known to be of type
D.}\label{tab:MainThm:espo:id}
\begin{tabular}{|p{4cm}|p{3cm}|p{3.5cm}|}

\hline  sporadic &  $t$ & {\bf Type of $\ell$ or\newline class name of $\oc_\ell^L$} \\

\hline  any &  $(t, |L|)= 1$,  $t$ odd & \textup{1A}  \\

\cline{2-3} &   2 & $\ord(\ell)=4$\\

\hline $T$, $J_2$, $Fi_{22}$, $Fi_{23}$, $Co_2$  &  odd & \textup{2A}  \\

\hline $B$ &  odd & \textup{2A, 2C}\\

\hline $Suz$ &  any & \textup{6B, 6C} \\

\hline
\end{tabular}
\end{center}
\end{table}


\section{Applications to the classification of pointed Hopf algebras}\label{sec:applications-pointed}

 We say that a finite group $G$ \emph{collapses} if for any finite-dimensional
    pointed Hopf algebra $H$, with $G(H) \simeq G$, necessarily $H\simeq\ku G$.
Some applications of the results on Nichols algebras presented here to the classification of Hopf algebras
need the following Lemma.

\begin{lema}\label{thlema:colapsar} \cite[Lemma 1.4]{AFGV} The following
 statements are equivalent:

 \begin{enumerate}
     \item\label{it:viejo1} If $0\neq V \in \ydg$, then $\dim \toba(V) = \infty$.
     \item\label{it:viejo2} If $V\in \ydg$ is \emph{irreducible},
         then $\dim \toba(V) = \infty$.
     \item\label{it:viejo3} $G$ collapses. \qed
 \end{enumerate}
\end{lema}

\begin{theorem}\label{teor:complete} \cite{AFGV, AFGV-espo}
Let $G$ be either an alternating group $\am$, $m\geq 5$, or a
sporadic simple group, different from
the Fischer group $Fi_{22}$, the Baby Monster $B$ and the Monster $M$.
Then $G$ collapses. \qed
\end{theorem}

The proof goes as follows: by the Lemma \ref{thlema:colapsar},
we need to show that $\dim \toba(V) = \infty$
for any irreducible $V = M(\oc, \rho)$.
If $\oc$ is of type D, this follows from Theorem \ref{th:racks-claseD};
and we know those classes of type D by Theorems \ref{th:racks-an-liquidados},
\ref{th:racks-liquidados}. The remaining
pairs $(\oc, \rho)$ are treated by abelian techniques,
namely one finds an abelian subrack, computes the
corresponding diagonal braiding arising from $\rho$ and applies \cite{H-all}.

\medbreak
However, there are finite non-abelian groups that do not collapse.
Furthermore, the classification of all finite-dimensional pointed
Hopf algebras with group $G$ is known, when $G=\s_3$, $\s_4$ or
$\mathbb D_{4t}$, $t\geq 3$, see \cite{AHS,GG,FG}, respectively.

\appendix

\section{ Examples of finite-dimensional Nichols
algebras}\label{subsec:nichols-examples-fd}

Table \ref{tab:app} contains several examples of pairs $(X, \q)$ such that
$\dim\toba(X,\q) < \infty$; we give the dimension, the top degree
and the reference where the example appeared\footnote{The Nichols
algebra corresponding to $\Q_{\Z/5, 2}$ was actually computed by
Mat\'\i as Gra\~na. The quadratic Nichols algebra corresponding to
$\Oc_{2}^5$ was computed by Jan-Erik Roos; Gra\~na showed that
this is a Nichols algebra. The computation of the Nichols algebras corresponding to $(\oc_{2}^n, \chi)$, $n =4,5$, was done in \cite{GG} using \textsf{Deriva} with the help of M. Gra\~na.}.

\begin{table}[ht]
\begin{center}
\caption{Finite-dimensional $\toba(X, \q)$}\label{tab:app}

\begin{tabular}{|c|c|c|p{2,5cm}|c|c|p{2cm}|}
\hline     $X$ & {\bf rk} & $\q$ & {\bf Relations} & $\dim
\toba(V)$ & {\bf top } & {\bf Ref. }

\\\hline $\D_3$ & 3 & -1  & 5 in degree 2     & $12 = 3.2^{2^{\quad}}$  & $4
= 2^2$ & \cite{MS}

\\\hline $\T$ & 4 & -1  & 8 in degree 2,\newline 1 in degree 6
& $72$  & $9 = 3^2$ & \cite{G-cm}

\\ \hline $\Q_{\Z/5, 2}$ & 5 &  -1 & 10  in degree 2,\newline  1 in degree 4
  & $1280 = 5.4^4$  & $16 = 4^2$ & \cite{AG-adv}

\\ \hline $\Q_{\Z/5, 3}$ & 5 &  -1 & 10  in degree 2,\newline  1 in degree 4
  & $1280 = 5.4^4$  & $16 = 4^2$ & dual of the\newline preceding

  \\ \hline $\Oc^4_2$ &  6& -1  & 16 in degree 2
  & $576 = 24^3$  & $12$ & \cite{FK, MS}

  \\ \hline $\Oc^4_2$ &  6& $\chi$  & 16 in degree 2
  & $576 = 24^3$  & $12$ & \cite{GG}

  \\ \hline $\Oc^4_4$ &  6& -1  & 16 in degree 2
  & $576 = 24^3$  & $12$ & \cite{AG-adv}

\\ \hline $\Q_{\Z/7, 3}$ &  7  & -1 & 21 in degree 2,\newline   1 in degree 6
  & $326592 = 7.6^6$  & $36 = 6^2$ & \cite{G-zoo}

\\ \hline $\Q_{\Z/7, 5}$ &  7  & -1 & 21 in degree 2,\newline   1 in degree 6
  & $326592 = 7.6^6$  & $36 = 6^2$ & dual of the\newline preceding

\\ \hline $\Oc^5_2$ & 10 & -1 & 45 in degree 2
  & $8294400$  & $40$ & \cite{FK, G-zoo}

\\ \hline $\Oc^5_2$ & 10 & $\chi$ & 45 in degree 2
  & $8294400$  & $40$ & \cite{GG}

\\ \hline
\end{tabular}

\end{center}
\end{table}


\section{Questions}\label{subsec:questions}

\begin{question-app}\label{que-app:rackssimples-typeD}
Determine all simple racks of type D.
\end{question-app}

\begin{question-app}\label{que-app:racks}
For any finite indecomposable rack $X$, for any $n\in \N$, and for
any $\q\in H^{2}(X, \mathbf{GL}(n,\ku))$, determine if $\dim \toba(X, \q)
< \infty$.
\end{question-app}

\begin{question-app}\label{que-app:collapses-degree1}
If $X$ collapses at $1$, does necessarily $X$ collapse?
\end{question-app}

\begin{question-app}\label{que-app:cociclo-diagonal-menosuno}
Compute all cocycles $\q\in Z^{2}(X,\ku^{\times})$ such that $\q_{ii} =  -1$.
\end{question-app}

\begin{question-app}\label{que-app:H2-trasp}
Is $H^2(\Oc^m_2, \ku^{\times}) \simeq \ku^\times\times\G_2$ for $m \geq 4$?
\end{question-app}

\begin{question-app}\label{que-app:simple-affine-subracks} Find the proper subracks of $\Q$.
\end{question-app}

\begin{question-app}\label{que-app:trasp-sm}
Are the cocycles $-1$ and $\chi$  equivalent by twist?
Recall that $H^{2}(\sm,\ku^{\times}) \simeq \Z/2$ for $m\geq4$, see \cite{schur}.
\end{question-app}

\medbreak\subsection*{Acknowledgements} N. A., G. A. G. and L. V.
want to thank the Organizing Committee the invitation to attend to
the \textit{XVIII Coloquio Latinoamericano de \'Algebra} and the
warm hospitality in S\~ao Pedro during the Colloquium. Specially,
G. A. G. thanks the invitation to give the mini-curse ``Quantum
Groups and Hopf Algebras''.

\end{document}